\documentclass[12pt]{amsart}
\usepackage[utf8]{inputenc}
\usepackage{setspace}
\usepackage[english]{babel}
\usepackage{mathabx}
\usepackage{framed}
\usepackage{soul}
\usepackage[all]{xy}
\usepackage[margin=1in]{geometry}
\usepackage{arydshln}
\usepackage[colorinlistoftodos]{todonotes}
\usepackage{pb-diagram}
\usepackage[mathscr]{euscript}
\usepackage{multicol}
\usepackage{stmaryrd} 
\usepackage{empheq} 
\usepackage{tikz-cd}
\tikzset{
  symbol/.style={
    draw=none,
    every to/.append style={
      edge node={node [sloped, allow upside down, auto=false]{$#1$}}}
  }
}

\usepackage{amsmath,amsfonts,amssymb,amsthm,mathrsfs,latexsym,mathtools}
\usepackage{commath}
\usepackage[T1]{fontenc}
\usepackage{hyperref}
\usepackage{caption} 
\usepackage{subcaption} 

\usepackage{lipsum}

\DeclareMathAlphabet{\mathbbmsl}{U}{bbm}{m}{sl}

\usetikzlibrary{matrix}

\newcommand{\C}{\mathbb{C}} 
\newcommand{\Z}{\mathbb{Z}}

\newcommand{\bbQ}{\mathbb{Q}}
\newcommand{\Hom}{\textup{Hom}}

\newcommand{\Gr}{\mathrm{gr}}

\newcommand{\Sym}{\textup{Sym}}
\newcommand{\SL}{\mathrm{SL}}
\newcommand{\NS}{\mathrm{NS}}

\newcommand{\bH}{\mathbb{H}}
\newcommand{\bP}{\mathbb P}
\newcommand{\bD}{\mathbb{D}}
\newcommand{\val}{\mathrm{val}}
\newcommand{\bN}{\mathbb{N}}

\usepackage{amscd} 

\usepackage{enumitem} 

\newtheorem{theorem}{Theorem}[section]
\newtheorem{definition}[theorem]{Definition}
\newtheorem{proposition}[theorem]{Proposition}
\newtheorem{corollary}[theorem]{Corollary}
\newtheorem{question}[theorem]{Question}

\newtheorem{lemma}[theorem]{Lemma}

\newtheorem{remark}[theorem]{Remark}

\newtheorem{thm}{Theorem}

\tikzset{commutative diagrams/.cd,
mysymbol/.style={start anchor=center,end anchor=center,draw=none}
}

\title{Elliptic-elliptic surfaces and the Hesse pencil}
\author{Fran\c{c}ois Greer and Yilong Zhang}
\dedicatory{Dedicated to Herb Clemens on the occasion of his 85th birthday}
\date{Sep 27, 2024}

\begin{document}

\maketitle

\begin{abstract}
    We construct a family of elliptic surfaces with $p_g=q=1$ that arise from base change of the Hesse pencil. We identify explicitly a component of the higher Noether-Lefschetz locus with positive Mordell-Weil rank, and a particular surface having maximal Picard number and defined over $\bbQ$. These examples satisfy the infinitesimal Torelli theorem, providing a second proof of the dominance of period map, which was first obtained in \cite{EGW}. A third proof is provided using the Shioda modular surface associated with $\Gamma_0(11)$. Finally, we find birational models for the degenerations at the boundary of the one-dimensional Noether-Lefschetz locus, and extend the period map at those limit points.
\end{abstract}

\section{Introduction}

An elliptic surface is a smooth projective surface $X$ with a map $f$ to a certain base curve $C$, whose general fiber is an elliptic curve. We always mark a distinguished section of the elliptic fibration map, which determines a commutative group law on the general fiber, and hence on the set sections; this is the Mordell-Weil group of the fibration $f$. We will construct certain special elliptic surfaces with $p_g=q=1$, whose Mordell-Weil groups and period images can be understood explicitly. These numerical invariants imply that the base curve has genus one, and the Hodge line bundle $L = f_*(\omega_f)$ has degree one. They are the simplest non-product elliptic surfaces of Kodaira dimension one, and they are parametrized by a 10-dimensional moduli space, $F_{1,1}$. Since $L$ has a unique nonzero section, up to scaling, there is a distinguished elliptic fiber of $f$ which represents the canonical class. For this reason, elliptic surfaces with $p_g=q=1$ are referred to as {\it elliptic-elliptic surfaces}.

In a sense, the most algebraic surfaces are those with maximal Picard number: 
$$\rho(X) = h^{1,1}(X).$$
All rational elliptic surfaces have the same (maximal) Picard number, but the most special ones are those which are {\it extremal}, meaning that they have torsion Mordell-Weil group, and Miranda and Persson \cite{MirPer86} classified them. Shioda and Inose \cite{ShiIno77} classified all K3 surfaces with maximal Picard number, referring to them as {\it singular}\footnote{The word ``singular'' here is used to mean ``very special''. All K3 surfaces that are singular in this sense are smooth and projective.}; by contrast, each such K3 surface admits an elliptic fibration structure with a section of infinite order. For surfaces of Kodaira dimension one, there are the Shioda modular surfaces \cite{Shi72}, which are universal families of elliptic curves with level structure. When the level subgroup of $\SL_2(\Z)$ is sufficiently small, the Shioda modular surface has Kodaira dimension one, but it always has torsion Mordell-Weil group. Persson constructed more examples \cite{Persson} with maximal Picard number by taking double cover of ruled surfaces, but again the Mordell-Weil groups are torsion. This brings us to our central question.




\begin{question}\label{question}
   Can one explicitly describe an elliptic surface of Kodaira dimension one with maximal Picard number and a section of infinite order?
\end{question}

From the point of view of Hodge theory, the Leray spectral sequence of an elliptic surface $f:X\to C$ degenerates at $E_2$ \cite{Zucker79}, and there is a decomposition of weight-two Hodge structures
\begin{equation}\label{eqn_H^2elliptic}
    H^2(X,\C)\cong H^0(C,R^2f_*\C)\oplus H^1(C,R^1f_*\C)\oplus H^2(C,f_*\C).
\end{equation}

The first and the last terms are Hodge-Tate, and they correspond to the class of zero section and the classes of curves supported on fibers, respectively. The middle term $H^1(C,R^1f_*\C)$ is the most interesting part. The integral $(1,1)$-classes of $H^1(C,R^1f_*\C)$ correspond to sections of $f:X\to C$ \cite{Zucker79,CoxZuc79}. When $X$ has Kodaira dimension one, there is a unique elliptic fibration structure, given by a pluricanonical morphism, so there is no ambiguity in referring to the Mordell-Weil of $X$ itself.

From general Noether-Lefschetz considerations, combined with dominance of the period map for $F_{1,1}$, we show in Prop. \ref{zariskidense} that there exist infinitely many elliptic-elliptic surfaces $X$ with any Mordell-Weil rank $r$ within the range $0\leq r \leq 10$, as allowed by Hodge theory. But the proof is not constructive, and it is often difficult to locate a section of infinite order.
We answer Question \ref{question} affirmatively with an explicit construction. 

\begin{theorem}\label{mainthm}
    There is an elliptic-elliptic surface $Y_p\to E$, defined over $\bbQ$, with two $I_3$ fibers and one $I_6$ fiber. Moreover, it has maximal Picard number $\rho=h^{1,1}=12$, its Mordell-Weil rank is 1, and it satisfies the infinitesimal Torelli theorem.
\end{theorem}

\begin{figure}[h]
\centering
\includegraphics[width=0.4\textwidth]{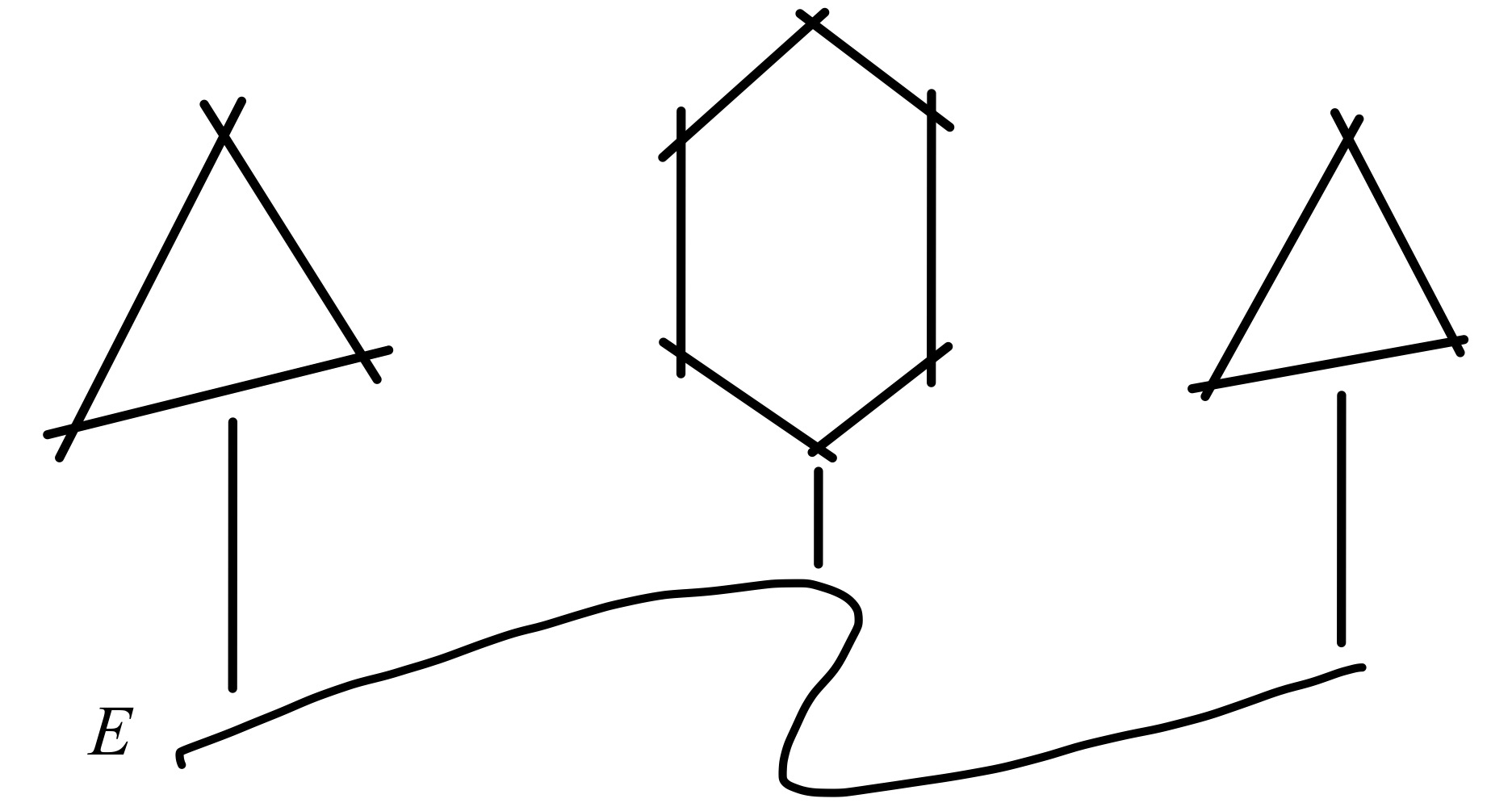}
\caption{The surface with maximal Picard number and fiber types $2I_3+I_6$}
\end{figure}

This is the first example of an elliptic-elliptic surface which satisfies the infinitesimal Torelli theorem. It differs from the examples studied in \cite{Ikeda}, which vary in one-parameter families with constant period point.
In \cite{EGW}, the period map for $p_g=q=1$ elliptic surfaces is shown to be dominant via a degeneration to the boundary of the moduli space $F_{1,1}$. The period map is locally injective in a neighborhood of the moduli point corresponding to the surface in Theorem \ref{mainthm}. This gives a second proof (see Corollary \ref{cor_IVHS}) of the dominance theorem via specialization to an interior point of $F_{1,1}$. We will also give a third proof (see Corollary \ref{shiodatorelli}) using a similar argument near the interior point of $F_{1,1}$ corresponding to the Shioda modular surface for $\Gamma_0(11)$.

\begin{corollary} (cf. \cite[Theorem 1.1]{EGW})
    The period map $P_{1,1}: F_{1,1}\to \bD/\Gamma$ from the moduli space of elliptic-elliptic surfaces is dominant.
\end{corollary}

Our construction begins by locating certain non-torsion multisections on the Hesse pencil that are invariant under an order 3 automorphism. Via base change, we produce a new elliptic surface over the normalization of a multisection, together with a section of infinite order. There is an induced $\Z_3$-action on this new elliptic surface, which descends to an isogeny on the base. The quotient produces interesting examples of elliptic-elliptic surfaces with large Picard rank and positive Mordell-Weil rank. For a particular multisection with tangency conditions at the special fibers, the construction produces an elliptic surface with maximal Picard rank.

The infinitesimal variation of the Hodge structure for an elliptic-elliptic surface $Y$ is closely related to the Kodaira-Spencer map at the canonical divisor $K_Y$, which can be read off from the resolution complex, and this allows us to deduce that the derivative of the period map is injective.


We locate the original trisections by expressing the Hesse pencil in terms of a relative Weierstrass equation and pulling back lines. We provide two interpretations of the family trisections; one is in terms of a birational involution on the $\Z_3$-quotient of the Hesse pencil (Section \ref{sec_involution}), and the other is in terms of Hesse configuration (Appendix \ref{sec_Hesse}).\\

\noindent\textbf{Noether-Lefschetz locus.}
Our second result concerns a higher Noether-Lefschetz locus in $F_{1,1}$ containing the example from Theorem \ref{mainthm}.
\begin{theorem}\label{thm_NL}
    There is a surface $S$ in the moduli space $F_{1,1}$ of elliptic-elliptic surfaces, a curve $T\subseteq S$ and a point $p\in T$, such that 
\begin{itemize}
    \item a general point of $S$ is an elliptic-elliptic surface with four singular fibers of type $I_3$, $\rho=10$, and Mordell-Weil rank $0$;
    \item a general point of $T$ is an elliptic-elliptic surface with four singular fibers of type $I_3$, $\rho=11$, and Mordell-Weil rank $1$;
     \item $p$ corresponds to the elliptic-elliptic surface from Theorem \ref{mainthm}, with two singular fibers of type $I_3$ and one of type $I_6$, $\rho=12$, and Mordell-Weil rank 1.
\end{itemize}
All the surfaces above satisfy the infinitesimal Torelli theorem.
\end{theorem}

These families arise from an alternative construction for the main example: we first pass to $\Z_3$-quotient of Hesse pencil, descending the trisection to the quotient surface $X'$, and then base change via a triple cover parametrized by a Hurwitz space. The subtlety here is that the base change involves semistable reduction of a type $IV^*$ singular fiber. For this reason, the base change does not increase the degree of the fundamental line bundle. \\



\noindent\textbf{Compactifying the period map.} We study to study the period image of the one-parameter family $T$ in Theorem \ref{thm_NL}, which is a modular curve by Noether-Lefschetz theory. There are three singular surfaces in the family, but only one of them produces a strictly mixed limiting Hodge structure of weight 2.


\begin{theorem}\label{thm_LMHS}
    The local monodromy actions for the one-parameter family $T$ at the three critical points $a=0,1,\infty$ have the following properties.
    \begin{itemize}
        \item At $a=0$, the monodromy actions on both $H^1$ and $H^2$ have finite order;
        \item At $a=1$, the monodromy action on $H^1$ (resp. $H^2$) has infinite (resp. finite) order;
         \item At $a=\infty$, the monodromy actions on both $H^1$ and $H^2$ have infinite order.
    \end{itemize}
    In particular, the period map associated to $H^2$ at $a=0,1$ is extendable to interior points of the period space.
\end{theorem}

We prove this by finding a normal crossing model of the family near the three critical points, and then use the Clemens-Schmid sequence to compute the limiting mixed Hodge structure.
Unlike for K3 surfaces, which are simply connected and whose degenerations are determined by the monodromy action on $H^2$ \cite{Kulikov77, Persson77}, for elliptic-elliptic surfaces, $H^1$ is nontrivial; the action on $H^1$ provides a further obstruction to the existence of a smooth semistable limit.

The case $a=0$ is analogous to a Type I degeneration of K3 surfaces. There is a smooth extremal elliptic-elliptic surface as the semistable limit, after taking a 3-to-1 base change, so the limiting mixed Hodge structure for both $H^1$ and $H^2$ are pure. At $a=1$, the semistable limit of the base elliptic curve is reducible of Kodaira type $I_3$, so the monodromy on $H^1$ has infinite order and obstructs a smooth semistable limit. However, the monodromy action on $H^2$ is still finite order — we find a smooth surface birational to a K3 as a component of a normal crossing limit, which forces the limiting mixed Hodge structure on $H^2$ to be pure. The case $a=\infty$ is analogous to a Type III degeneration of K3 surfaces. We find a semistable limit whose components are rational surfaces, with dual complex homotopy equivalent to a once-pinched $2$-torus, i.e., a Kodaira fiber of type $I_1$.\\

\noindent\textbf{Structure of the paper.}
In Section \ref{sec_prelim}, we recall the definition of Mordell-Weil group for an elliptic surface, Shioda-Tate formula, and explain the behaviors of sections and singular fibers under base change. In the end, we discuss basic facts about the Hesse pencil. 
In Section \ref{sec_RES}, we discuss the quotient of the Hesse pencil by a certain automorphism of order 3. We explain how to obtain an elliptic-elliptic surface with four $I_3$ fibers by taking certain base change of the quotient.
In Section \ref{sec_onepara}, we construct a one-parameter family of non-torsion trisections of the Hesse pencil which are invariant under the $\Z_3$-action. This produces the one-parameter family $T$ of elliptic-elliptic surfaces with positive Mordell-Weil rank. We find an elliptic surface with maximal Picard rank in the family.
In Section \ref{sec_properties}, we discuss some general properties of elliptic-elliptic surfaces, including the canonical fiber, torsion sections, primitivity of the non-torsion section, $J$-invariants of the base elliptic curves $E_a$, and fields of definition.
In Section \ref{sec_period}, we introduce the moduli space of elliptic-elliptic surfaces and its associated period domain of weight 2. We show that there are infinitely many elliptic-elliptic surfaces with a given Mordell-Weil rank, and define higher Noether-Lefcshetz loci.
In Section \ref{sec_IVHS}, we prove the infinitesimal Torelli theorem for the elliptic-elliptic surfaces we constructed in Section \ref{sec_RES}. This gives a proof of Theorem \ref{thm_NL}. 
In Section \ref{hurwitz}, we use a certain double Hurwitz space to formally construct the families appearing in Theorem \ref{thm_NL} and complete the proof.
In Section \ref{sec_Shioda}, we find the canonical fiber of the Shioda modular surface for $\Gamma_0(11)$ and show that this surface also satisfies the infinitesimal Torelli theorem.
In Section \ref{sec_involution}, we study a birational involution of the rational elliptic surface $X'$ that interchanges trisections and fibers, and resolve the indeterminacy explicitly on a blow up.
In Section \ref{sec_compactify}, we study a compactification of the one-parameter family $T$, and compute the limit mixed Hodge structure at the three critical values. This proves Theorem \ref{thm_LMHS}.\\




\noindent\textbf{Acknowledgements}. The authors are grateful to Donu Arapura for introducing them to the problem and for many helpful discussions. The authors would like to thank Dori Bejleri, Lian Duan, Philip Engel, Anatoly Libgober, Xiyuan Wang, and Jaroslaw Wlodarczyk for useful communication. The first author was partially supported by NSF grant DMS-2302548.

\tableofcontents

\section{Preliminaries}\label{sec_prelim}
\subsection{Shioda-Tate formula}

For an elliptic surface $f:X\to C$ with a chosen section, the {\it Picard number} $\rho(X)$ denotes the rank of the Neron-Severi group
$$\NS(X) = \mathrm{span} \{[D]\in H^2(X,\Z)|D\subset X\ \textup{is a curve}\}.$$
Its Mordell-Weil group is the group of sections of $f$, relative to the chosen zero section, or equivalently the group of rational points of the elliptic curve defined over the function field of $C$. The Mordell-Weil rank  is the rank of the Mordell-Weil group. 

The Shioda-Tate formula \cite[Corollary VII.2.4]{Miranda} relates the Picard number $\rho(X)$ and the Mordell-Weil rank $r$ of an elliptic surface:
\begin{equation}\label{eqn_ShiodaTate}
    \rho(X)=2+r+\sum_{t\in C}(n_t-1),
\end{equation}
where $n_t$ is the number of irreducible components in the fiber $f^{-1}(t)$. Since the generic fiber is irreducible, this is a finite sum. 

In short, the formula says that the Neron-Severi group is generated by classes of curves supported in fibers, and classes of section curves. 
The Shioda-Tate formula \eqref{eqn_ShiodaTate} is also a consequence of the Leray decomposition \eqref{eqn_H^2elliptic}. We refer to \cite{Ara,AraSol} for related discussion.


\subsection{Non-torsion multisections}
Let $f:X\to C$ be an elliptic surface with zero section $C_0$. Then there is a degree $n$ morphism
$$[n]: X\to X$$
whose restriction to a smooth fiber $f^{-1}(t)=E_t$ is the multiplication by $n$: $E_t\to E_t$.

\begin{definition}
   A multisection is an irreducible curve $D\subset X$ such that the degree of $f|_D:D\to C$ is $\geq 2$.
\end{definition}
\begin{definition}
    A section or multisection is called torsion if it is an irreducible component of the preimage $T_n:=[n]^{-1}(C_0)$ for some $n\in\bN$.
\end{definition}


\begin{proposition}\label{prop_ntmultisection}
    If $D\subset X$ is a multisection, and moreover $f|_D$ is ramified along a smooth fiber of $f$, then $D$ is not torsion.
\end{proposition}
\begin{proof}
    Let $f_U:X_U\to U$ be the restriction of $f$ to smooth fibers. Then $T_n$ is the closure of $T_n\cap X_U$. However, $T_n\cap X_U$ is \'etale over $U$ of degree $n^2$, so any torsion multisection must have branch points contained in $C\smallsetminus U$.
\end{proof}

Observe that the property of being torsion or non-torsion is preserved under base change. As a consequence of the Shioda-Tate formula \eqref{eqn_ShiodaTate}, the Mordell-Weil rank of an elliptic surface increases upon finite base change if and only there is a non-torsion multisection  that becomes a non-torsion section. This will be used in Section \ref{sec_onepara}.


\subsection{Base change and reduction of singular fibers}
The possible singular fiber types were classified by Kodaira, N\'{e}ron, and Tate. We refer to \cite[Section 5.4]{SchShi19} for the list of singular fibers. In what follows, we will be mostly focused on $I_n$ fibers, a loop consisting of $n$ rational curve components, and the $IV^*$ fiber, a nonreduced fiber with 7 components arranged in an $\tilde{E}_6$ configuration.

The base change of an elliptic surface ramified along a singular fiber gives a singular surface. The minimal desingularization typically results in a different Kodaira fiber type, since one needs to perform birational modifications on the total space to obtain a new smooth model. The correspondence is given in \cite[p.105, Table 2]{SchShi19}. As an example, the $d$-to-1 base change of an $I_n$ fiber produces an $I_{dn}$ fiber:
\begin{equation}\label{eqn_I_nbasechange}
     I_{dn}\xrightarrow{d:1} I_n.
\end{equation}

In particular, a 2-to-1 base change branched along an $I_3$ fiber creates creates a surface with three $A_1$ singularities. Resolving these singularities, we obtain a smooth model with an $I_6$ fiber, as pictured below.

\begin{figure}[h]
\centering
\includegraphics[width=0.5\textwidth]{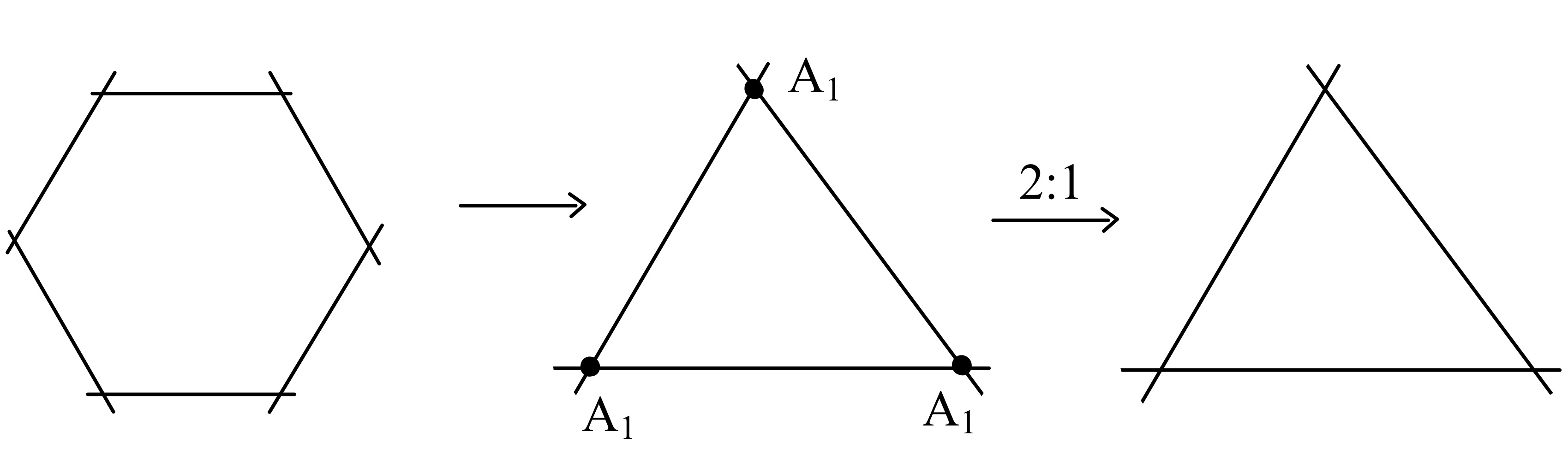}
\caption{quadratic base change of $I_3$}
\end{figure}

Unlike for $I_n$ fibers, which are semistable, base change along a non-reduced fiber results in a surface with non-normal total space. One must normalize and perform blow-ups and blow-downs to get a relatively minimal regular model. For example, the 2-to-1 and 3-to-1 base change ramified along an $IV^*$ fiber are given below.
\begin{align}\label{eqn_IV*basechange}
   IV\xrightarrow{2:1} IV^*,\  I_0\xrightarrow{3:1} IV^*
\end{align}

This process is closely related to semistable reduction for curves of arithmetic genus 1. We will explicitly compute the semistable reduction for the base change of a type $IV^*$ fiber in the Appendix \ref{sec_app}.

\subsection{Hesse pencil}
The Hesse pencil $f:X\to \mathbb P^1$ is the universal family of elliptic curves with level three structure. It is the unique rational elliptic surface with four $I_3$ fibers \cite{Bea82}. Its Mordell-Weil group is isomorphic $\Z_3\oplus \Z_3$, so it is extremal. We refer to \cite{BarHul85} for a projective model of the Hesse pencil.



The automorphism group $G$ of the Hesse pencil has order of $216$, and fits into an exact sequence \cite[Proposition 4.1]{ArtDol}:
$$0\to \ker(\phi)\to G\xrightarrow{\phi}\textup{Aut}(\mathbb P^1).$$

The image of $\phi:G\to Aut(\bP^1)$ is isomorphic to the alternating group $A_4$, realized as the even permutations of the four critical values in $\bP^1$. The kernel of $\phi$ is a group of order 18 generated by the Mordell Weil group $\Z_3\oplus \Z_3$ and the fiberwise involution $x\mapsto -x$.
We are interested in a particular subgroup of order 3 in $G$, whose image under $\phi$ cyclically permutes three singular fibers.

By contracting two of the ($-2$)-curves in each singular fiber, Hesse pencil admits a Weierstrass model given by equation \cite[p.167, \# 68]{SchShi19}
\begin{equation}\label{eqn_number68}
   y^2+(3tx+t^3-1)y=x^3. 
\end{equation}

The coordinate on the base $\bP^1$ is $t$, and the fibration is induced by coordinate projection. This is an affine and singular model of the Hesse pencil; we use the transformation \eqref{eqn_number68-s} to view the fiber at $\infty$. The singular fibers are irreducible nodal cubics, but the total space has $A_2$ singularities located at the nodes in the singular fibers. It has singular fibers over $1,\zeta_3,\zeta_3^2,\infty$, where $\zeta_3= e^{2\pi i/3}$.

\begin{figure}[h]
\centering
\includegraphics[width=0.5\textwidth]{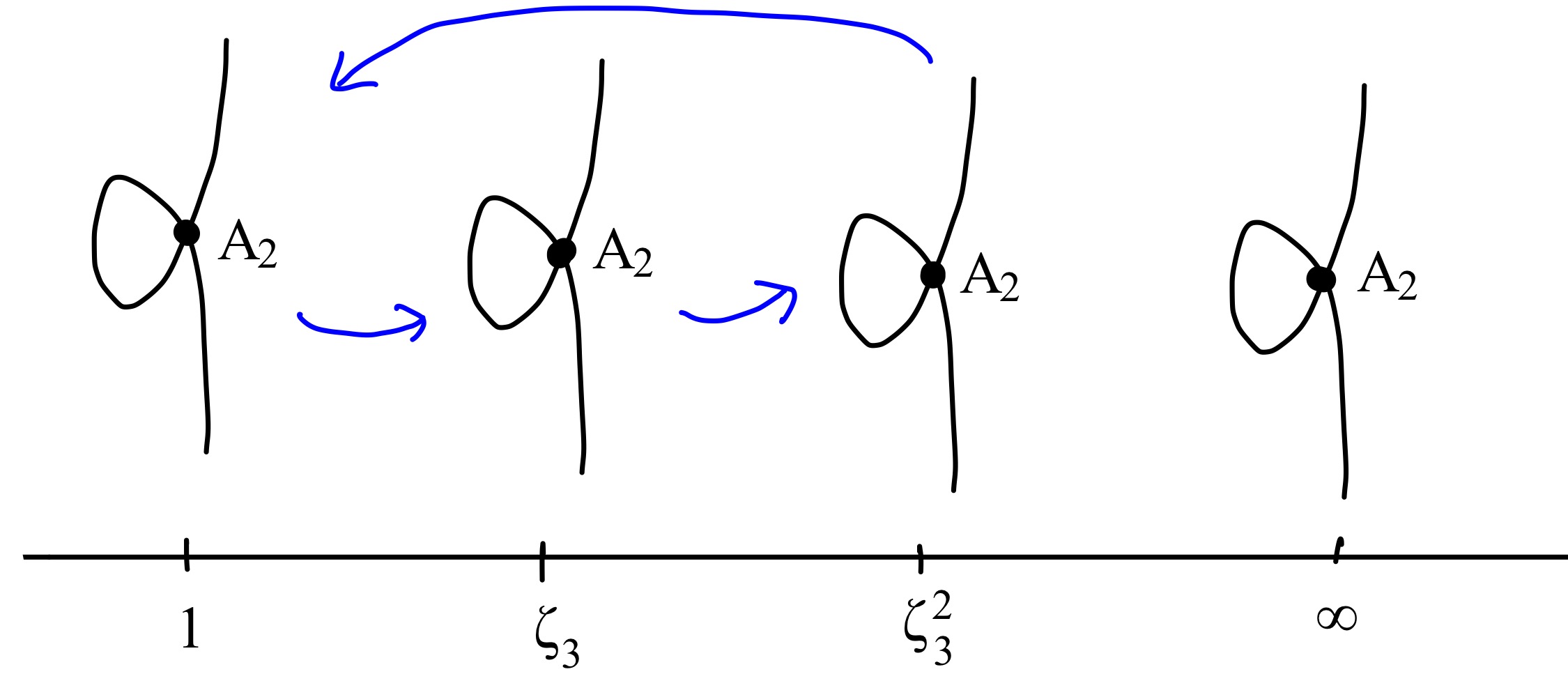}
\caption{Weierstrass model of the Hesse Pencil and $\Z_3$-action}
\end{figure}

Our order 3 subgroup of $G$ is generated by the automorphism
\begin{equation}\label{eqn_Z3-t}
    \sigma: (x,y,t)\mapsto (\zeta_3^{-1}x,y,\zeta_3t).
\end{equation}
It permutes the fibers over $1,\zeta_3,\zeta_3^2$ and fixes the fiber over $\infty$, and it will preserve certain non-torsion trisections of the Hesse pencil.




\section{Elliptic-Elliptic Surfaces with Four $I_3$ Fibers}\label{sec_RES}

In this section, we construct elliptic-elliptic surfaces with four singular fibers of type $I_3$ by taking a quotient of the Hesse pencil by an automorphism of order 3, and then performing a certain base change and birational modification.

\subsection{$\Z_3$-quotient of Hesse pencil}
Let $X\to \mathbb P^1$ be the Hesse pencil. According to the formula \eqref{eqn_Z3-t}, the automorphism $\sigma$ preserves the elliptic fibration structure, so the relatively minimal regular model $X'$ of the quotient $X/\Z_3$ is again an elliptic surface over $\bP^1$.

\begin{lemma}   \label{lemma_X'}
$X'\to \mathbb P^1$ is an extremal rational elliptic surface with three singular fibers of Kodaira types $IV^*,I_3,I_1$. More precisely, there is a commutative diagram:

\begin{equation}\label{diagram_RES}
\begin{tikzcd}
I_0, \ 3\times I_3,I_3 \arrow[r, symbol=\subset] \arrow[d] &X\arrow[d,dashrightarrow, "/\Z_3"] \arrow[r] &  \mathbb P^1 \arrow[d,"/\Z_3"]\\
IV^*, \ I_3,\ I_1 \arrow[r, symbol=\subset]& X'\arrow[r]  & \mathbb P^1.
\end{tikzcd}
\end{equation}

\end{lemma}
\begin{proof}
   We still work with the Weierstrass equations. The order three automorphism \eqref{eqn_Z3-t} identifies the three singular fibers on the $t$-chart. In particular, it identifies the three $I_3$ fibers on the smooth model.
   
   The action stabilizes the $t=0$ fiber and has three fixed points there. $f^{-1}(0)$ is the elliptic curve $y^2-y=x^3$, which has an automorphism:  $x\mapsto \zeta_3x$. The three fixed points are given by $(x,y)=(0,0),(0,1)$, and the point at infinity, and the quotient is $\mathbb P^1$ with multiplicity three, and the quotient surface $X/\Z_3$ has three isolated singularities of type $A_2$ located along this nonreduced fiber. The minimal desingularization gives the desired $IV^*$ fiber. 
    

 To study the fiber at infinity, we apply the transformation  
 \begin{equation}\label{eqn_transform_s}
  x'=s^{2}x,\ y'=s^3y,\  \textup{and}\ s=1/t
 \end{equation}
to the Weierstrass equation \eqref{eqn_number68}, so on $s$-chart, the Weierstrass equation becomes
 \begin{equation}\label{eqn_number68-s}
     (y')^2+3x'y'+(1-s^3)y'-(x')^3=0.
 \end{equation}

It has an $A_2$-singularity over $s=0$. However, the quotient will be smooth there; the $\Z_3$-action \eqref{eqn_Z3-t} on the $s$-chart becomes $\sigma: (x',y',s)\mapsto (x',y',\zeta_3^{-1}s)$, so via the substitution $v=s^3$, the $vy'$ term appears and implies that the nodal point of the $v=0$ fiber located at $x'=-1,y'=1$ is not a singularity on the surface. Hence, the smooth model $X'$ of the quotient surface has an $I_1$ fiber over $\infty$.
\end{proof}


\begin{lemma}\label{lemma_MWX'}
    The Mordell-Weil group of $X'$ is isomorphic to $\Z_3$.
\end{lemma}

  This can be read off from the table in \cite[p.166, row 69]{SchShi19}. Explicitly, the Hesse pencil has nine torsion sections. Under the $\Z_3$-action, three of them are fixed and descend to sections of $X'\to \mathbb P^1$. The remaining six sections form two disjoint $\Z_3$-orbits. Their images are two torsion trisections on $X'$ totally branched at $0$ and $\infty$.

\begin{remark}\normalfont\label{rmk_P1-to-P1}
 Conversely, one can view the Hesse pencil as the relative smooth model of the base change of $X'\to \mathbb P^1$ with respect to 3-to-1 map $\mathbb P^1\to \mathbb P^1$ totally branched at $0$ and $\infty$. This point of view is taken in \cite[Theorem 7.3]{MirPer86}.
\end{remark}

\subsection{Base change of $X'$} Let $E\to \mathbb P^1$ be the triple cover totally branched at $0,\infty$, and simply branched at $r_1,r_2\in \C^*$, with $\C^*=\mathbb P^1\smallsetminus\{0,\infty\}$. By the Riemann-Hurwitz formula, $E$ is a genus one curve.



\begin{proposition}\label{prop_basechangeX'}
   Let $Y$ be the relatively minimal regular model of the base change $X'\times_{\mathbb P^1}E$. Then $Y\to E$ is an elliptic-elliptic surface. For general $r_1,r_2$, $Y$ has four singular fibers of type $I_3$, and $\rho(Y)\geq 10$.
\end{proposition}
\begin{proof}
According to Lemma \ref{lemma_X'}, the rational elliptic surface $X'$ has singular fibers of types $IV^*$ and $I_1$ at $0$ and $\infty$. Since $E\to \mathbb P^1$ triply ramified at these two points, their relatively minimal regular models are $I_0$ and $I_3$ respectively (cf. \eqref{eqn_I_nbasechange} and \eqref{eqn_IV*basechange}). For a general choice of simple branch points, $r_1,r_2$ are disjoint from the remaining critical point below the $I_3$ fiber of $X'$, so its pullback consists of three copies of $I_3$. Finally, the degree of fundamental line bundle is $1$ by computing the Eular characteristic $\chi_{top}(Y)=4\times \chi_{top}(I_3)=12$, which implies the elliptic surface has $p_g=q=1$ (cf. \cite[p.36]{Miranda}). By the Shioda-Tate formula, $\rho(Y)=2+r+4\times 2\ge 10$.
\end{proof}
This construction gives a two-parameter family $S$ of elliptic-elliptic surfaces, varying $r_1,r_2$. In Section \ref{hurwitz}, we will show that the base of the family is a double Hurwitz space, and it dominates a surface in $F_{1,1}$ as well as in the period space, using the results of Section \ref{sec_IVHS}.

\subsection{Commuting base change and quotient}
Let $\tilde{E}$ be the normalization of the fiber product $E\times_{\mathbb P^1}\mathbb P^1$, with respect to the 3-to-1 base change $\mathbb P^1\to \mathbb P^1$ totally branched at $0$ and $\infty$ (cf. Remark \ref{rmk_P1-to-P1}), then $\tilde{E}\to E$ is unramified, and therefore an isogeny between elliptic curves, while the other projection $\tilde{E}\to \mathbb P^1$ is 3-to-1 and simply branched at the 6 pre-images of $r_1$ and $r_2$. Consider the base change $X\times_{\mathbb P^1}\tilde{E}$ of Hesse pencil, with the group $\Z_3$ acting diagonally. Let $\tilde{X}\to \tilde{E}$ be the relatively minimal regular model.


\begin{proposition}
    The $\Z_3$-action on $X$ extends to a free action on the elliptic surface $\tilde{X}\to \tilde{E}$. The quotient is birational to the elliptic-elliptic surface $Y\to E$ (cf. Proposition \ref{prop_basechangeX'}). 
\end{proposition}
\begin{proof}
    For general $r_1,r_2$, the $\tilde{E}\to \mathbb P^1$ has a branch point lying under smooth fibers of $X$. In that case, the fiber product $X\times_{\mathbb P^1}\tilde{E}$ is smooth. (Otherwise, if one of $r_1,r_2$ lies under the $I_3$ fiber, one needs to resolve $A_1$ singularities on the fiber product.) Since the $\Z_3$-action on $\tilde{E}$ is free, its action on the fiber product is free. The quotient agrees with $Y\to E$ on smooth fibers. Hence by uniqueness, its relatively minimal regular model is exactly the elliptic-elliptic surface $Y$.
\end{proof}

As a consequence, there is a commutative diagram of rational maps between elliptic surfaces, regular over the loci in the base curves where the fibers are smooth.
 \begin{equation}\label{diagram}
\begin{tikzcd}
\tilde{X}/\tilde{E}\arrow[d,"/\Z_3",dashed] \arrow[r,dashed] &  X/\mathbb P^1 \arrow[d,"/\Z_3",dashed]\\
Y/E\arrow[r,dashed]  & X'/\mathbb P^1.
\end{tikzcd}
\end{equation}

This point of view will be useful in the next section, where we construct elliptic surfaces with sections of infinite order, and in Section \ref{sec_properties} for finding the canonical fiber of $Y$.

\section{Elliptic-Elliptic Surfaces with Positive Mordell-Weil Rank}\label{sec_onepara}

In this section, we construct a one-parameter subfamily $T$ of elliptic-elliptic surfaces with generic Mordell-Weil rank equal to one. It will be a component of the Noether-Lefschetz locus for the two-parameter family from the previous section.

\subsection{Trisections of the Hesse pencil}
Finding sections of infinite order on $Y$ requires finding certain multisections on the Hesse pencil. Our goal is to show:


\begin{proposition}\label{prop_Da}
    There exists a family of irreducible curves $\{D_a\}$ for $a\in \C^*\smallsetminus \{1\}$ on the Hesse pencil such that 
    \begin{enumerate}
    \item $D_a\to \mathbb P^1$ has degree 3;
        \item $D_a$ is $\Z_3$-invariant, and the quotient $E_a:=D_a/\Z_3$ is a genus one curve;    
        \item $D_a$ is non-torsion;
        \item For $a=4$, $D_4$ is tangent to the singular fibers over $t=1,\zeta_3,\zeta_3^2$. 
    \end{enumerate}
\end{proposition}



The proof will be divided into a series of propositions. The idea is that the Weierstrass equation \eqref{eqn_number68} may be regarded as a cubic pencil of cubic curves. In other words, \eqref{eqn_number68} defines a hypersurface in $\C^3$, whose projection to the $x,y$ coordinate has degree three. Since a line in $\C^2$ intersects a general cubic curve at three points, it corresponds to a multisection of degree three over the base $\bP^1$, or a \textit{trisection}. 

\begin{definition}\normalfont\label{def_Da}
    Let $D_a$ be the curve on the Hesse pencil defined by the equation $y=a$ in the Weierstrass equation \eqref{eqn_number68}.
\end{definition}
\begin{proposition}\label{prop_trisection}
    The curve $\{D_a\}$ is a non-torsion trisection on the Hesse pencil $X$, and it has geometric genus one. Moreover, the composition of the desingularization $\tilde{D}_a\to D_a$ and $f$
    is simply branched at 6 points which lie in two $\Z_3$-orbits. When $a=4$, one such orbit lies over $1,\zeta_3,\zeta_3^2\in \bP^1$.

\end{proposition}

\begin{proof}
   Plugging $y=a$ into the Weierstrass equation \eqref{eqn_number68} yields an affine equation for the trisection $D_a$:
    \begin{equation}\label{eqn_y=a}
       x^3+A(a,t)x+B(a,t)=0, 
    \end{equation}
    where $A(a,t)=-3at$ and $B(a,t)=-at^3+a-a^2$. The discriminant of the cubic \eqref{eqn_y=a} is given by $4A(a,t)^3+27B(a,t)^2=0$. Assuming $a\neq 0$ and choosing a branch for $\sqrt{a}$, this simplifies to
\begin{equation}\label{eqn_y=a-disc}
    (t^3-(\sqrt{a}-1)^2)(t^3-(\sqrt{a}+1)^2)=0.
\end{equation}
As long as $a\neq 0,1$, this equation has six simple roots that are permuted by the $\Z_3$-action $t\mapsto \zeta_3 t$. Hence, $D_a$ has 6 simple ramification points in two $\Z_3$-orbits. The affine curve \eqref{eqn_y=a} is smooth as long as $a\neq 0,1$. Observe that it has a $t^3-1$ factor iff $a=4$, so $D_a$ is tangent to the three $I_3$ fibers at $t^3=1$ iff $a=4$.
Equation \eqref{eqn_y=a} defines a smooth affine cubic (when $a\in \C^*\smallsetminus\{1\}$), so its smooth compactification $\tilde{D}_a$ is a smooth genus one curve. By the Riemann-Hurwitz formula, there are no further branch points, so $\tilde{D}_a\to \mathbb P^1$ is unbranched at $\infty$. By Proposition \ref{prop_ntmultisection}, $D_a$ is non-torsion.
\end{proof}


\begin{proposition}\label{prop_Sa_zerosection}
Let $D_a$ be the trisection above, with $a\in \C^*\smallsetminus \{1\}$.
\begin{enumerate}
    \item $D_a$ has an ordinary triple point over $t=\infty$ and no other singularities, so it has arithmetic genus $g_a(D_a)=4$.
    \item $D_a$ intersects each singular fiber of $f:X \to \bP^1$ away from singularities of the fiber, in the same irreducible component as the zero section.
    \item $D_a$ is disjoint from the zero section.

\end{enumerate}
\end{proposition}
\begin{proof}
Applying the coordinate change \eqref{eqn_transform_s} to \eqref{eqn_y=a}, we find that the trisection has equations 
\begin{equation}\label{eqn_S-s}
   y'=as^3,\ \textup{and}\  (a^2-a)s^6+3ax's^3+as^3=(x')^3
\end{equation}
in the affine chart near $t=\infty$. Its homogeneous degree 3 part is $as^3-(x')^3=0$, so as long as $a\neq 0$, $D_a$ has an ordinary triple point at $s=0$ and has three local branches. It has Milnor number $\mu=4$, then by Milnor-Jung, $\delta=(4+3-1)/2=3$, and the arithmetic genus $g_a$ is the geometric genus plus $\delta$, which is $4$.

    For (2), both \eqref{eqn_y=a} and \eqref{eqn_S-s} have degree 3 in $x$ and $x'$ meaning that the multisection is disjoint from the zero section, which is placed at infinity. The four singular fibers of \eqref{eqn_number68-s} are at $s=0,1,\zeta_3,\zeta_3^2$ with singularities at $(s,x',y')=(0,-1,1), (1,0,0), (\zeta_3,0,0),(\zeta_3^2,0,0)$. One can verify directly that \eqref{eqn_S-s} does not pass through these points. By the construction of the relatively minimal regular model from the Weierstrass equation, their intersection with the singular fiber lies in the same component with the zero section.
\end{proof}


\begin{proposition}\label{prop_trisection/Z3}
    The trisection $D_a$ for $\C^*\smallsetminus \{1\}$ is invariant under the automorphism $\sigma$. Moreover, the induced action $\Z_3$ on $\tilde{D}_a$ is free, and the quotient $E_a:=\tilde{D}_a/\Z_3$ is a smooth curve of genus one. In particular, $E_a$ descends to a trisection on the rational elliptic surface $X'$ (cf. Lemma \ref{lemma_X'}). The 3-to-1 map $E_a\to \mathbb P^1$ is totally branched at $0$ and $\infty$ and simply branched at two other points.
\end{proposition}
\begin{proof}
   The $\Z_3$-action $\sigma|_D: (x,4,t)\mapsto (\zeta_3^{-1}x,4,\zeta_3t)$ preserves the trisection, and is fixed point free when $t\neq \infty$. Near $t=\infty$, $\sigma$ permutes the three branches, so the induced action $\sigma:\tilde{D}_a\to \tilde{D}_a$ is free, and $\tilde{D}_a\to \tilde{D}_a/\Z_3$ is an isogeny of elliptic curves.    

   Over $t=0$, the three points on $\tilde{D}_a$ are in a single $\Z_3$-orbit, so $E_a$ is totally branched at $t=0$. Over $t=\infty$, the $\Z_3$-action fixes the ordinary triple point, and its image in $E_a$ is the only point over $\infty$, so $E_a$ is totally ramified there. Lastly, by the proof of Proposition \ref{prop_trisection}, there are two simple branch points at $(1\pm \sqrt{a})^2$.
\end{proof}

\begin{remark}\normalfont
    One should note that our construction of trisections by pulling back a line depends on the choice of the Weierstrass equation. If we transform the Weierstrass equation from \eqref{eqn_number68} to \eqref{eqn_number68-s} and then projects to the $x,y$ coordinate, the image of $D_a$ is a cubic curve, instead of a line (it intersects the blowup center, which is a length 9 subscheme, with multiplicity 6, so its lift to the blow up still has degree three). Pulling back the invariant lines to \eqref{eqn_number68-s} and then taking $\Z_3$-quotients of the resulting trisections gives rational curves instead of genus one curves.
\end{remark}

\subsection{A section of infinite order} We now apply the base change construction of Section \ref{sec_RES}, using the genus one curves $E_a \to \bP^1$ that are trisections of $X'$.
\begin{proposition}\label{prop_rho=11family}
    The elliptic-elliptic surface $Y_a\to E_a$ arising from each trisection $E_a$, $a\in\C^*\smallsetminus\{1\}$, has Mordell-Weil rank at least one, and $\rho(Y_a)\geq 11$. 
\end{proposition}

\begin{proof}
    The base change $X\times_{\mathbb P^1}\tilde{D}_a\to \tilde{D}_a$ has a tautological section given by the diagonal map. Its birational image in the relatively minimal regular model $\tilde{X}\to \tilde{D}_a$ extends to a non-torsion section. The $\Z_3$-action on $\tilde{X}$ preserves the section, and so it descends to a non-torsion section of $Y_a\to E_a$.
    
   Alternatively, one can first take the $\Z_3$-quotient and then base change to construct the elliptic surface $Y_a\to E_a$, as in the diagram \eqref{diagram}. By Proposition \ref{prop_trisection/Z3}, $E_a=\tilde{D}_a/\Z_3$ is a genus one curve, and the map $E_a\to \mathbb P^1$ is 3-to-1 totally branched at $0$ and $\infty$, and simply branched at at two other points $r_1,r_2=(1\pm \sqrt{a})^2$. The base change  $X'\times_{\mathbb P^1}E_a$ has a section of infinite order, and by Proposition \ref{prop_basechangeX'}, its relatively minimal regular model is an elliptic-elliptic surface. By the Shioda-Tate formula, 
   $$\rho(Y_a)=2+r+4\times 2= 10+r\geq 11,$$
   since $r\ge 1$.
\end{proof}





\subsection{Elliptic surfaces with maximal Picard rank}
The following proves Theorem \ref{mainthm}.

\begin{proposition}\label{prop_maximalPicard1}
    There exists a particular surface $Y_p$ in the one-parameter family $\{Y_a\}$ with maximal Picard number and Mordell-Weil rank one. It has three singular fibers of types $$2\times I_3+I_6.$$
\end{proposition}
\begin{proof}
Take the trisection $D_4$ in Proposition \ref{prop_Da}. Its $\Z_3$-quotient $E_4\subset X'$ is simply branched at $t=1$ where the fiber of $X' \to \bP^1$ is $I_3$, so after base change, the total space acquires has three new $A_1$ singularities. After resolving these singularities, cf. \eqref{eqn_I_nbasechange}, the elliptic surface $Y_4$ has an $I_6$ fiber and two $I_3$ fibers. By the Shioda-Tate formula, 
$$\rho(Y_4)=2+r+2\times 2+5=11+r\le h^{1,1}(Y_4)=12.$$ 
On the other hand, $Y_4 \to E_4$ has a non-torsion section, so $r=1$ and $\rho(Y_4)=12$. 
\end{proof}

In addition, one can identify an extremal elliptic surface, which arises as a semistable limit of the one-parameter family as $a\to 0$ (cf. Proposition \ref{prop_LMHS_0}).

\begin{proposition}\label{prop_maximalPicard2}
    There exists an elliptic-elliptic surface $Y_0$ with maximal Picard number and Mordell-Weil rank equal to zero, with two singular fibers of type 
    $$I_3+I_9.$$
\end{proposition}
\begin{proof}
    Take the base change of $X'\to \mathbb P^1$ with respect to the Galois cover $\mathbb P^1\to \mathbb P^1$ totally branched at $0,1,\infty$. In the relatively minimal regular model $Y_0$, the $IV^*$ fiber becomes a smooth fiber, the $I_3$ becomes an $I_9$, and the $I_1$ becomes an $I_3$. By the Shioda-Tate formula, 
    $$\rho(Y_0)=2+r+8+2=12+r\le h^{1,1}(Y_0)=12,$$
    so $r=0$ and $\rho(Y_0)=12$.
\end{proof}

Both of these examples are defined over $\bbQ$ (cf. Proposition \ref{prop_numberfield}). By \cite[Lemma 4.4]{Ara19}, the Tate conjecture holds for these surfaces.

\section{Properties}\label{sec_properties}
In this section, we discuss some properties of the special elliptic-elliptic surfaces constructed above. First, all elliptic surfaces $f:Y\to E$ with $p_g=q=1$ share the following:
\begin{enumerate}
\item The fundamental line bundle $f_*(\omega_f)$ has degree one, and its pullback is the canonical bundle of $Y$ \cite[Theorem 5.44]{SchShi19}. 
\item The topological Eular characteristic $\chi_{top}(Y)=12$.
\item The morphism $f:Y\to E$ induces an isomorphism on fundamental groups $\pi_1(Y)\to \pi_1(E)\cong \Z^2$ \cite[Chapter II, Proposition 2.1]{FM}.
\item The morphism $f:Y\to E$ induces an isomorphism of weight one Hodge structures $H^1(E,\Z)\to H^1(Y,\Z)$ \cite[Lemma VII.1.1]{Miranda}.
\item Since $h^{2,0}(Y)=1$, there is a unique fiber of $f:Y\to E$ which is the zero locus of a nontrivial holomorphic two form; we call this the \textit{canonical fiber}.
\end{enumerate}

From now, we restrict to those elliptic-elliptic surfaces $Y\to E$ which occur in the two-parameter family of Proposition \ref{prop_basechangeX'}. Recall that the $J$-map for an elliptic surface associates to $t\in E$ the $J$-invariant of the fiber $Y_t$ in $\bP^1$.
\begin{proposition}
   The degree of the $J$-map $E\to \bP^1$ associated to $Y\to E$ is equal to 12.
\end{proposition}
\begin{proof}
  This is true for all elliptic surfaces with semistable fibers and fundamental line bundle of degree one. Alternatively, the degree of $J$-map is 12 for the Hesse pencil $X$, 4 for the $\Z_3$-quotient surface $X'$, and then 12 again for $Y$ because of the 3-to-1 base change.
\end{proof}

\subsection{Canonical fiber} The canonical classes of the surfaces $Y\to E$ constructed in Proposition \ref{prop_basechangeX'} can be identified in terms of the 3-to-1 map $E\to \bP^1$ used for base change.


\begin{proposition}\label{prop_p+r-infty}
    The canonical class $K_Y$ is represented by the pullback of the divisor 
    $$\tilde{r}_1+\tilde{r}_2-\tilde{\infty}$$
    from $E$, where $\tilde{r}_1$ and $\tilde{r}_2$ are the two simple ramification points on $E$ lying over $r_1$ and $r_2$, and $\tilde{\infty}$ is the total ramification point on $E$ lying over $\infty\in \bP^1$.
\end{proposition}
\begin{proof}
We use the notation of \eqref{diagram}. The base change $X\times_{\mathbb P^1}\tilde{E}$ of the Hesse pencil is normal and has only $A_n$ singularities, so passing to the crepant resolution $\tilde{X}$ does not affect the canonical class. We first compute the canonical class $K_{\tilde{X}}$, which is the pullback of $K_X$ plus ramification divisors. Since $K_X$ is represented by $-F_t$ for any point $t\in \mathbb P^1$ we can take $t=\infty$. The preimage of $\infty\in \bP^1$ in $\tilde{E}$ consists of three points, denoted $\tilde{\infty}_1,\tilde{\infty}_2,\tilde{\infty}_3$, which form a single $\Z_3$-orbit. The six simple ramification points of $\tilde{E}\to \mathbb P^1$ fall into two $\Z_3$-orbits over $r_1$ and $r_2$. Denote them by $\tilde{r}_{ij}$, with $i=1,2$, $j=1,2,3$.
    

Therefore, $K_{\tilde{X}}$ is represented by pullback of the divisor $\sum_{i,j}\tilde{r}_{ij}-\tilde{\infty}_1-\tilde{\infty}_2-\tilde{\infty}_3$ on $\tilde{E}$. The support of this divisor falls into three $\Z_3$-orbits, and they correspond to $\tilde{r}_1,\tilde{r}_2,\tilde{\infty}$ on the $\Z_3$-quotient curve $E$. Finally, since the $\Z_3$-action on $\tilde{X}$ is free, $K_Y$ is represented by the pullback of $\tilde{r}_1+\tilde{r}_2-\tilde{\infty}$ on $E$.
\end{proof}

\begin{remark}\normalfont
    Alternatively, this can be proved via semistable reduction (cf. Appendix \ref{sec_app}).



\end{remark}




\begin{corollary}\label{cor_K}
The canonical fiber of any $Y$ arising from the base change construction of Proposition \ref{prop_basechangeX'} is the fiber over $\tilde{0}\in E$, the point of total ramification lying over $0\in \bP^1$.
\end{corollary}
\begin{proof}
    We apply Riemann-Hurwitz to the triple cover $g:E \to \bP^1$ with simple ramification at $\tilde{r}_1,\tilde{r}_2\in E$ and total ramification at $\tilde{0},\tilde{\infty}\in E$. Using $K_{\bP^1}=-0-\infty$, we obtain
    \begin{align*}
        0 &= K_E = g^*K_{\bP^1}+\tilde{r}_1+\tilde{r}+2\cdot\tilde{0}+2\cdot\tilde{\infty}\\
        &= \tilde{r}_1+\tilde{r}_2 - \tilde{0} - \tilde{\infty}.
    \end{align*}
The result now follows from Proposition \ref{prop_p+r-infty}.
\end{proof}

In particular, the canonical fiber $K_Y$ is always isomorphic to the Eisenstein elliptic curve with $J$-invariant zero, for all $Y$ in the two-parameter family of elliptic-elliptic surfaces.

\subsection{Torsion Mordell-Weil groups}
\begin{proposition}\label{prop_MWtor}
Any elliptic-elliptic surface $Y$ that arises from Proposition \ref{prop_basechangeX'} has 
$$MW(Y)_{tor}\cong \Z_3.$$
\end{proposition}
\begin{proof}
   The Mordell-Weil group of the rational elliptic surface $X'$ is isomorphic to $\Z_3$ (cf. Lemma \ref{lemma_MWX'}), by the classification of extremal surfaces. Pulling back sections from $X'$, we obtain $\Z_3\subset MW(Y)$. It remains to show there are no other torsion sections. If $\tau\subset Y$ were such a torsion section, then its image under the map $Y \to X'$ would be a torsion trisection $\tau'\subset X'$. However, the triple cover $\tau' \to \bP^1$ must have the same branching data as $E \to \bP^1$. Since $(r_1,r_2)\neq (1,1)$ the trisection $\tau'$ must be ramified along a smooth fiber of $X' \to \bP^1$, so by Proposition \ref{prop_ntmultisection}, it is non-torsion. 
\end{proof}   

\begin{remark}
    Even if we include the limit case $Y_0$ of Proposition \ref{prop_maximalPicard2}, effectively allowing $(r_1,r_2)=(1,1)$, Proposition \ref{prop_MWtor} remains true. From the singular fiber types, the Neron-Severi lattice of $Y_0$ has discriminant $3\times 9/|MW_{tor}(Y)|^2$ \cite[Corollary 6.39]{SchShi19}. Since $|MW_{tor}(Y)|\ge 3$, the discriminant must be equal to $3$, and $MW_{tor}(Y)$ cannot be any larger.
\end{remark}
   
   
   

\subsection{Primitivity of the tautological section}
From here to the end of this section, we restrict further to the one-parameter subfamily of elliptic-elliptic surfaces with positive Mordell Weil rank constructed in Proposition \ref{prop_rho=11family}.

\begin{corollary} 
A general member $Y_a$ of the one-parameter family $T$ has
$$MW(Y_a)\cong \Z\oplus \Z_3.$$
\end{corollary}

One might ask whether the non-torsion section given by the diagonal map $E_a\to X'\times_{\mathbb P^1}E_a$ generates the free part of the Mordell-Weil group.


\begin{proposition}
    The tautological section of $Y_a \to E_a$ induced by trisection $E_a\subset X'$ is primitive in the group $MW(Y_a)$.
\end{proposition}
\begin{proof}
    Let $\Sigma$ be the non-torsion section of $Y_a\to E_a$ constructed in Proposition \ref{prop_rho=11family}, which passes through the same reducible fiber components as the zero section, $C_0$. The powers of any section with respect to $*$, the Mordell-Weil group law, give sections whose classes can be computed explicitly. Assume that $\Sigma = \Sigma_0^{*n}$ for some $n>1$.  First, we have
    \begin{equation}\label{mwgplaw}
        [\Sigma_0^{*n}] = n\Sigma_0 - (n-1)C_0 + kF + R,
    \end{equation}
    where $R$ is a divisor supported on the singular fiber components of $Y_a \to E_a$ which do not meet $C_0$. First, we have
$$
        0=\Sigma \cdot C_0 = n \Sigma_0\cdot C_0 + n-1+k,
$$
    so $k = -n(h+1)+1$, where $h = \Sigma_0\cdot C_0\geq 0$. On the other hand,
    $$-1 = \Sigma\cdot \Sigma = -n^2-(n-1)^2+R\cdot R + 2n\Sigma_0\cdot R - 2n(n-1)h + 2k;$$
    \begin{equation}\label{theeqn}
        2n^2 (h+1)-2 = R\cdot R + 2n\Sigma_0\cdot R.
    \end{equation}
Now, $R$ is a linear combination of $(-2)$-curve classes orthogonal to $C_0$ and $F$, which generate a negative definite sublattice of $\NS(Y_a)$. If $\Sigma_0\cdot R=0$, then the RHS of equation \ref{theeqn} is nonpositive, producing a contradiction. If $\Sigma_0\cdot R\neq 0$, then one checks using \ref{mwgplaw} that $3|n$ and $3|(R\cdot R+2n\Sigma_0\cdot R)$. This gives a contradiction by reducing equation \ref{theeqn} modulo 3.
\end{proof}

Therefore, the Mordell-Weil lattice for $Y_a$ in the sense of \cite{SchShi19} is $A_1$.

\subsection{$J$-Invariant of the base elliptic curve}
Here, we compute the $J$-invariant of the base elliptic curves $E_a$ for $Y_a$ in the one-parameter family. By substitution $z=s^3$ in \eqref{eqn_S-s}, we find that $E_a$ has affine equation $(a^2-a)z^2+3ax'z+az=(x')^3$, which is a smooth elliptic curve when $a\neq 0, 1,\infty$. A direct computation shows that 
\begin{equation}\label{eqn_JEa}
   J(E_a)=\frac{27 a (a + 8)^3}{(a -1 )^3}. 
\end{equation}

In particular, it defines a degree-four map $\C^*\smallsetminus\{1\}\to M_{1,1}$ to the coarse moduli space of elliptic curves. The elliptic surface with maximal Picard number in Theorem \ref{mainthm} has base corresponding to $a=4$ with $J$-invariant $4\cdot 1728$. The extremal elliptic surface in Proposition \ref{prop_maximalPicard2} has base curve with $J$-invariant 0. In fact, it corresponds to the limit as $a\to 0$ (cf. Corollary \ref{cor_afamily}), and of course $\lim_{a\to 0}J(E_a)=0$.


\subsection{Field of definition} In the end, we restrict to the two elliptic-elliptic surfaces $Y_p$ and $Y_0$ with maximal Picard rank from Proposition \ref{prop_maximalPicard1} and \ref{prop_maximalPicard2}.
\begin{proposition}\label{prop_numberfield}
  The elliptic surfaces $Y_p$ and $Y_0$ are defined over $\bbQ$.
\end{proposition}
\begin{proof}
 We first use the description of $Y$ in Proposition \ref{prop_basechangeX'} by base change of the rational elliptic surface $X'$ and then perform birational modifications. Since both $X'\to \mathbb P^1$ and $E_4\to \mathbb P^1$ are defined over $\bbQ$, their fiber product is defined over $\bbQ$.

  In an affine chart, we look at the tensor product of rings
  \begin{equation} \label{eqn_Y_p_overQ}
      \C[x',y',z]/(f)\otimes_{\C[z]} \C[w,z]/(g)\cong \C[x',y',z,w]/(f,g),
  \end{equation}
  where $f=(y')^2+3x'y'+(1-z)y'-x^3$ and $g=w^3-12zw-4z-12z^2$. It has two $A_2$ singularity at $(-1,1,0,0)$ and at $(0,0,1,4)$, and an $A_5$ singularity at $(0,0,1,-2)$. Since they are all $\bbQ$-points, the minimal resolution of \eqref{eqn_Y_p_overQ} has a $\bbQ$ structure and is an open subspace of $Y_p$ away from canonical fiber $F_{\tilde{0}}$. 
  
  To find the other chart covering $Y_p$, we use the equation \eqref{eqn_y=a} in $t$-chart, and find that the base change $D_4\times_{\mathbb P^1}X$ of Hesse pencil has $\bbQ$-structure and defined by 
  \begin{equation} \label{eqn_Y_p_overQ2}
      \C[x,y,t]/(f_1)\otimes_{\C[z]} \C[w,t]/(f_2)\cong \C[x,y,z,w]/(f_1,f_2),
  \end{equation}
  where $f_1=y^2+(3tx+t^3-1)y-x^3$, $f_2=12+12tw+4t^3-w^3$.

Then the $\Z_3$ invariant of subring is generated by $x^3,y,z^3,w^3,xz,wz,x^2w,xw^2$, and $f_1$ and $f_2$ in terms of these monomials. The relations of these minomials have $\bbQ$-coefficient, therefore, the resolution of the spectrum of the $\Z_3$ invariant subring of \eqref{eqn_Y_p_overQ2} provides an open subspace of $Y_p$ containing the canonical fiber and is defined over $\bbQ$.

The argument for $Y_0$ is similar.
\end{proof}




\section{Moduli and Period Spaces}\label{sec_period}

\subsection{Moduli spaces} 
For a given genus $g\geq 0$ and degree $d>0$, one can consider the moduli stack $F_{g,d}$ of elliptic surfaces $f:X \to C$ with $g=g(C)$ and $d = \deg f_*(\omega_f)$. Each surface $X$ has a Weierstrass model $w(X)$ obtained by contracting all components of the fibers which do not meet the marked section. The fibers of $w(X)\to C$ are irreducible, but this surface can have ADE singularities. The moduli stack $W_{g,d}$ of Weierstrass models has the advantage of being expressible as a global quotient. We assume for convenience that $2d+1>g$. Let $\mathcal P_{g,d}$ be the universal Picard scheme of degree $d$ over $\mathcal M_g$, which has dimension $4g-3$. Let $\mathcal L_{g,d}$ be the universal line bundle on $\mathcal P_{g,d}\times_{\mathcal M_g} \mathcal C_g$. Consider the pushforward
$$V=p_{1*}(\mathcal L_{g,d}^{\otimes 4}\oplus \mathcal L_{g,d}^{\otimes 6})$$
on $\mathcal P_{g,d}$ which parametrizes Weierstrass coefficients $(A,B)\in H^0(C,L^{\otimes 4}\oplus L^{\otimes 6})$. By Serre duality, and the fact that $d$ is sufficiently large, $H^1(C,L^{\otimes 4}\oplus L^{\otimes 6})=0$, so $V$ is locally free of rank $10d-2g+2$ by the Theorem on Cohomology and Base Change. Let $U\subset V$ be the open subset of pairs $(A,B)$ such that the surface
$y^2 = x^3+Ax+B$
has at worst ADE singularities. By \cite{Kas}, this is equivalent to the condition that $4A^3+27B^2 \not\equiv 0$, and for all $x\in C$,
$$\min(3\,\val_x(A), 2\,\val_x(B)) <12,$$
an open condition. Now, $\C^\times$ acts on $V$ with weights $(2,3)$ with respect to the direct sum decomposition, and 
$$W_{g,d} \cong U/ \C^\times,$$
an open substack of a weighted projective bundle over $\mathcal P_{g,d}$. By counting the dimensions above, it is a Deligne-Mumford stack of dimension $10d+2g-2$.

The stacks $F_{g,d}$ and $W_{g,d}$ have isomorphic sets of $\C$-points, since an elliptic surface $X$ is the minimal resolution of its Weierstrass model $w(X)$. As observed by Artin \cite{artin}, the stack $F_{g,d}$ is not separated, but there is a natural isomorphism $F_{g,d} \to W_{g,d}$ inducing the isomorphism on $\C$-points. Under this isomorphism, elliptic surfaces $X$ with a reducible fiber correspond to Weierstrass models $w(X)$ that have singularities. This locus forms a divisor $\Delta$ in moduli, which will be the pre-image under the period map of a certain Noether-Lefschetz divisor; see below for the definition of the weight 2 period map.

\subsection{Period maps}
For the remainder of this section, we will focus primarily on the case where $(g,d)=(1,1)$. These surfaces are called elliptic-elliptic. For each $Y\in F_{1,1}$, its middle cohomology lattice is isometric to $I_{3,11}$, the odd unimodular lattice of signature $(3,11)$. Since $H^{1,1}(Y,\Z)$ always contains a rank 2 lattice generated by the fiber and zero section classes $\{F,C_0\}$, it is standard to consider the {\it primitive cohomology} lattice
$$H^2_{\mathrm{prim}}(Y,\Z)= \langle F,C_0 \rangle ^\perp \subset H^2(Y,\Z),$$
which is isometric to $I\!I_{2,10}$, and supports a weight 2 polarized Hodge structure of type $(1,10,1)$. The space of such Hodge structures on $I\!I_{2,10}$ is a Hermitian symmetric domain $\bD = \bD(I\!I_{2,10})$ of type IV:
\begin{equation}\label{perioddomain}
\bD(I\!I_{2,10}) := \{ \C\omega\in \bP(I\!I_{2,10}\otimes \C): \omega\cdot \omega = 0,\, \omega\cdot \overline{\omega}>0 \}
\end{equation}
The arithmetic quotient $\bD/\Gamma$, where $\Gamma$ is the group of integral isometries $I\!I_{2,10} \to I\!I_{2,10}$. has the structure of a quasi-projective variety by \cite{bb}.

Any family $\mathcal Y \to B$ of elliptic-elliptic surfaces over a base scheme $B$ has a classifying map $B \to F_{1,1}$. For any point $b\in B$, the fiber $\mathcal X_b$ has a polarized Hodge structure on its primitive middle cohomology lattice, giving a weight 2 PVHS on $B$. Choosing an isometry $H^2_{\mathrm{prim}}(Y_b,\Z)\cong I\!I_{2,10}$ allows us to associate to $b$ a point in $\bD$; any two choices of isometry differ by an element of $\Gamma$. If $\widetilde{B}$ denotes the universal cover of $B^{an}$, we obtain a marked period map $\widetilde{B} \to \bD$ that is holomorphic by \cite{griffiths}. This marked period map fits into a commutative square
$$\xymatrix{
\widetilde{B} \ar[r]\ar[d] & \bD \ar[d] \\
B \ar[r] & \bD/\Gamma
}$$
The property above allows us to define a holomorphic period map $P_{1,1}$
$$P_{1,1}: F_{1,1} \to \bD/\Gamma$$
Since both sides are also algebraic, and $P_{1,1}$ is a definable morphism, the period map is algebraic. The main result of \cite{EGW} shows that $P_{1,1}$ is generically finite and dominant. The following questions are natural in comparison with the analogous situation for $(g,d)=(0,2)$, the case of elliptic K3 surfaces.
\begin{question}\label{quest1}
    What is the image of $P_{1,1}$ in $\bD/\Gamma$? What is the degree of this period mapping?
\end{question}
\begin{question}
    Can the ramification locus of $P_{1,1}$ be described via a geometric condition on the elliptic-elliptic surfaces $X\in F_{1,1}$?
\end{question}
\begin{question}
    Is the period map $P_{1,1}$ defined over $\bbQ$?
\end{question}
The constructions in this paper provide the first explicit examples of points in $F_{1,1}$ where $P_{1,1}$ is unramified, even though this must be true generically. As $\dim(F_{1,1}) = \dim(\bD/\Gamma)=10$, such points give independent proofs of Thm 1.1 in \cite{EGW}.

\begin{remark}\normalfont
    Apart from the trivial case $(g,d)=(0,1)$, where the period space is a point, the only other case of elliptic surfaces where the period space is quasi-projective is $(g,d)=(0,2)$. There, the answers to Questions 3.1-3.3 are well-known: $P_{0,2}$ is an isomorphism by the Global Torelli Theorem for lattice-polarized K3 surfaces, and $P_{0,2}$ is defined over $\bbQ$.
\end{remark}

\subsection{Noether-Lefschetz Loci}
The period space $\bD/\Gamma$ for elliptic-elliptic surfaces contains numerous special subvarieties which are normalized by arithmetic quotients of smaller type IV Hermitian symmetric domains. These are the (higher) Noether-Lefschetz loci. 
\begin{definition}
\normalfont    Let $\Lambda$ be an even integral lattice of signature $(2,n)$. We define the associated period domain as in \ref{perioddomain}:
    $$\bD(\Lambda) := \{\C \omega\in \bP(\Lambda\otimes \C): \omega\cdot \omega=0,\, \omega\cdot \overline{\omega}>0\}$$
    Given a primitive lattice embedding $\iota:\Lambda\hookrightarrow I\!I_{2,10}$, we have an associated inclusion $\bD(\Lambda) \hookrightarrow \bD(I\!I_{2,10})$. If $\Gamma(\Lambda)$ denotes the group of integral isometries $\Lambda \to \Lambda$, then $i$ induces
    $$\iota_*: \bD(\Lambda)/\Gamma(\Lambda) \to \bD/\Gamma,$$
    since every isometry of $\Lambda$ extends to an isometry of $I\!I_{2,10}$. Any image of $\iota_*$ for some $(\Lambda, \iota)$ is called a higher Noether-Lefschetz locus.
\end{definition}
In the extreme case when $n=0$, we obtain the CM points of $\bD/\Gamma$. If $d\in \bN$ is the discriminant of a rank 2, positive definite lattice with an embedding into $I\!I_{2,10}$, then the associated point of $\bD$ has coordinates in the imaginary quadratic field $\bbQ(\sqrt{-d})$. Even for fixed $d$, such points are dense in $\bD$ with respect to the classical topology. Conversely, 

\begin{proposition}\label{zariskidense}
    There exist infinitely many elliptic-elliptic surfaces $X$ with Mordell-Weil rank $r$ for any integer $0\leq r \leq 10$, such that all the fibers of $f:X \to C$ are irreducible.
\end{proposition}
\begin{proof}
    Since the CM points are dense in $\bD$, there images are also dense in the quotient $\bD/\Gamma$ with respect to the Zariski topology. The period map $P_{1,1}: F_{1,1} \to \bD/\Gamma$ is dominant, so the image is a Zariski open subset. The same is true after removing the divisor $\Delta$ parametrizing surfaces with a reducible fiber. Hence, there exist infinitely many elliptic-elliptic surfaces $X$ of maximal Picard rank such that all the fibers of $f:X \to C$ are irreducible. By the Shioda-Tate formula, the Mordell-Weil rank is equal to 10 for these surfaces. Fix one such surface $X_0$, and choose a lift $p\in \bD$ of its period point. This corresponds to a lattice $\Lambda_0$ of signature $(2,0)$ with a primitive embedding into $I\!I_{2,10}$. For each value of $s$ with $0\leq s \leq 10$, there exists an intermediate lattice
    $$\Lambda_0\subset \Lambda_s \subset I\!I_{2,10},$$
    by taking $\Lambda_s = \Lambda_{s-1}+w$ where $w\in \Lambda_{s-1}^\perp$ is primitive with $(w,w)<0$. If $N$ is a sufficiently small neighborhood of $p$, then $\bD(\Lambda_s)\cap N$ maps to a locally closed analytic subspace of $\bD/\Gamma$ with dimension $s$ which contains $P_{1,1}([X_0])$ and is disjoint from the period image of $\Delta$. A generic surface $X_s$ with period point in this subspace has Mordell-Weil rank $10-s$, by the Shioda-Tate formula again.
\end{proof}
\begin{question}
    What are the possible torsion subgroups of the Mordell-Weil groups for elliptic-elliptic surfaces?
\end{question}

Cox gave a complete answer to the analogous question in the cases $(g,d)=(0,1)$ and $(g,d)=(0,2)$ of rational elliptic surfaces and elliptic K3 surfaces, respectively, after determining all the possible ranks  \cite{coxtorsion}. The proof relies heavily on the Global Torelli Theorem for K3 surfaces; compare with Question \ref{quest1}.

The main example of this paper features nested families of elliptic-elliptic surfaces over $p\in T\subset S$ whose periods dominate higher Noether-Lefschetz loci in $\bD/\Gamma$, associated to nested lattices $\Lambda_0\subset \Lambda_1\subset \Lambda_2$. To determine these particular lattices of signature $(2,s)$, it suffices to find the transcendental lattices for the generic member of the nested families.

\begin{proposition}
    The transcendental lattice $T_Y$ of a generic surface in the family $S \to F_{1,1}$ (resp. $T$, $p\to F_{1,1}$) has the following invariants.
    \begin{itemize}
        \item For a general $Y\in S$, $T_Y\cong \Lambda_2$ has signature $(2,2)$ and discriminant group $\Z_3\oplus \Z_3$.
        \item For a general $Y\in T$, $T_Y\cong\Lambda_1$ has signature $(2,1)$ and discriminant group $\Z_3\oplus \Z_6$. 
        \item $T_{Y_p}\cong \Lambda_0$ has signature $(2,0)$ and discriminant group $\Z_2\oplus \Z_6$.
    \end{itemize}
\end{proposition}
\begin{proof}
    Let $Y$ be an elliptic-elliptic surface. The Neron-Severi lattice $N_Y$ inside $H^2(Y,\Z)$ has orthogonal complement $T_Y$, so the two discriminant groups are isomorphic. Using \cite{SchShi19}, $N_Y$ is an index $|MW(Y)_{tor}|$ overlattice of $I\!I_{1,1}\oplus A_2^{\oplus 4}$ (resp. $I\!I_{1,1}\oplus A_2^{\oplus 4}\oplus A_1$, $I\!I_{1,1}\oplus A_2^{\oplus 2}\oplus A_5\oplus A_1$) in the cases we consider here. Passing to a finite index overlattice changes the discriminant group by passing to $K^\perp/K$, where $K$ is an isotropic order 3 subgroup of the discriminant. The result now follows from the fact that $A_n$ has discriminant group $\Z/(n+1)\Z$.
\end{proof}

\section{Infinitesimal Variation of Hodge Structure}\label{sec_IVHS}

In this section, we study the derivative of period map $P_{1,1}$. In other words, we consider the infinitesimal variation of the weight 2 Hodge structures on $H^2(Y,\Z)$ for the elliptic-elliptic surfaces constructed in Proposition \ref{prop_basechangeX'} as they vary in the entire moduli space $F_{1,1}$. Following \cite[Proposition 5.1]{Ikeda}, we want to study the map
$$\phi: H^1(Y,T_Y)\to \Hom(H^{0}(\Omega_Y^2)\to H^1(\Omega^1_Y)).$$
The infinitesimal Torelli theorem in this context is equivalent to the injectivity of $\phi$.
\begin{lemma} \label{lemma_IVHS} 
Let $Y\to E$ be an elliptic-elliptic surface. If the canonical divisor $K$ of $Y$ is a smooth fiber, then $\phi$ is either injective or its kernel is 1-dimensional. The latter case occurs iff the Kodaira-Spencer map
$$\kappa: H^0(K,N_{K/Y})\to H^1(K,T_K)$$
is the zero. In other words, infinitesimal Torelli holds iff $\kappa$ is an isomorphism.
\end{lemma}
\begin{proof}
First, $\phi$ dualizes to 
$$\mu: H^{0}(\Omega_Y^2)\otimes H^1(\Omega^1_Y)\to H^1(\Omega^2_Y\otimes \Omega^1_Y)$$
The kernel of $\phi$ is dual to the cokernel of $\mu$, so infinitesimal Torelli holds iff $\mu$ is surjective.
Since $H^0(\Omega_Y^2)$ is one-dimensional, the unique section $s$ vanishes at the canonical fiber $K$. We have the exact sequence 
$$0\to \Omega_Y^1\xrightarrow{s} \Omega^2_Y\otimes \Omega^1_Y\to (\Omega^2_Y\otimes \Omega^1_Y)|_K\to 0$$
of sheaves, which gives the following long exact sequence on cohomology groups:
$$H^1(\Omega_Y^1)\xrightarrow{s} H^1(\Omega^2_Y\otimes \Omega^1_Y)\to H^1(K,(\Omega^2_Y\otimes \Omega^1_Y)|_K)\to H^2(\Omega_Y^1)\to H^2(\Omega^2_Y\otimes \Omega^1_Y)\to 0.$$
The image of $\mu$ coincides with the image of $s$. Since a non-isotrivial elliptic surface has no infinitesimal automoprhisms \cite[Theorem D]{Infiaut}, $H^0(T_Y)=0$, and by Serre duality $H^2(\Omega^2_Y\otimes \Omega^1_Y)\cong H^0(T_Y)^{\vee}$. Note also that $\dim H^2(\Omega_Y^1)=h^{1,2}(Y)=1$, so from the long exact sequence, infinitesimal Torelli holds iff $\dim H^1(K,(\Omega^2_Y\otimes \Omega^1_Y)|_K)=1.$
Using the assumption that $K$ is smooth with the adjunction formula, $\Omega^2_Y|_K\cong \Omega^1_K\cong \mathcal{O}_K$, so it suffices to consider $H^1(K,(\Omega^1_Y)|_K)$. The conormal sequence for $K\subset Y$ reads
$$0\to N_{K/Y}^{\vee}\to (\Omega^1_Y)|_K\to \Omega_K^1\to 0,$$
and it gives rise to the following long exact sequence:
\begin{equation}\label{eqn_exact}
    H^0(K,\Omega_K^1)\xrightarrow{\alpha} H^1(K,N_{K/Y}^{\vee})\to  H^1(K,(\Omega^1_Y)|_K)\to H^1(K,\Omega_K^1)\to 0
\end{equation}
Obverse that $\alpha$ is dual to the Kodaira-Spencer map at the canonical fiber as an elliptic curve:
$$\kappa: H^0(K,N_{K|Y})\to H^1(K,T_K).$$
Since both sides are 1-dimensional, $\kappa$ is either zero or an isomorphism. By the exactness of \eqref{eqn_exact} and the fact that $H^1(K,\Omega^1_K)\cong \C$, we conclude that $\kappa$ is an isomorphism iff $\dim  H^1(K,(\Omega^1_Y)|_K)=1$ iff infinitesimal Torelli holds for the surface $Y$.
\end{proof}

\begin{corollary}\label{cor_IVHS}
    Infinitesimal Torelli holds for all the elliptic-elliptic surfaces $Y$ that were constructed in Proposition \ref{prop_basechangeX'}.
\end{corollary}
    \begin{proof}
Corollary \ref{cor_K} implies that the canonical fiber $K_Y$ is smooth with $J$-invariant zero. By Lemma \ref{lemma_IVHS}, it suffices to show the Kodaira-Spencer map $\kappa$ at the canonical fiber is nonzero.
     Recall from the construction via base change and then $\Z_3$-quotient that $\tilde{E}\to \mathbb P^1$ is \'etale at $t=0$. The base change produces three copies of an analytic disk containing that fiber of the Hesse pencil, and the $\Z_3$-quotient identifies these three copies, so a neighborhood of the canonical fiber in $Y$ is analytically isomorphic to a neighborhood of the corresponding fiber in the Hesse pencil. Since the Hesse pencil is extremal, by Lemma \ref{lemma_KSextremal} below, the Kodaira-Spencer map is nonzero at all smooth fibers.  
     
     Alternatively, one can show that the $J$-map of the family of elliptic curves $Y\to E$ is triply ramified at the canonical fiber. Since the $J$-map locally factors through the stack $\mathcal{M}_{1,1}\to M_{1,1}$, and the latter map is triply ramified at $J=0$, the map $E\to \mathcal{M}_{1,1}$ is unramified at canonical fiber, so $\kappa$ is nonzero there.
    \end{proof}
    The lemma below is a special case of \cite{Ara,Ara19}. We provide a self-contained proof.
\begin{lemma}\label{lemma_KSextremal}
    Let $f: X\to C$ be a non-isotrivial extremal elliptic surface. Assume that all fibers are semistable (type $I_n$), and let $D\subset C$ be the set of critical values. Then the map
    \begin{equation}\label{eqn_KSfamily}
        \bar{\nabla}: f_*\omega_{X/C}\to R^1f_*\mathcal{O}_{X}\otimes \Omega^1_{C}(\log D)
    \end{equation}
    given by cup product and contraction of Kodaira-Spencer class, is an isomorphism.
\end{lemma}
\begin{proof}
    Extremality implies $MW(X/C)$ has rank zero, which by \cite[Lemma 2.2, 2.6]{Ara} says
    \begin{equation}\label{eqn_GrF1}
        \Gr_F^1H^1(C,R^1f_*\C)=0.
    \end{equation}
   By the degeneration of the Hodge-to-de Rham spectral sequence at $E_1$ \cite{Zucker79}, we have 
   $$\Gr_F^1H^1(C,R^1f_*\C)=H^1(C,\Gr_F^1DR(\mathcal{V})),$$
   where $\mathcal{V}:=R^1f_*\C$ and $DR(\mathcal{V})$ is the de Rham sequence
    $$\mathcal{V}\otimes \mathcal{O}\xrightarrow{\nabla} \mathcal{V}\otimes \Omega^1_{C}(\log D).$$
Its first graded piece is exactly \eqref{eqn_KSfamily}. Since the elliptic surface is non-isotrivial,  \eqref{eqn_KSfamily} is nonzero, and therefore an injective morphism of locally free sheaves, so 
$$\Gr_F^1DR(\mathcal{V})\cong \textup{coker}(\bar{\nabla})[-1].$$
Therefore, $H^1(C,\Gr_F^1DR(\mathcal{V}))=H^1(C,\textup{coker}(\bar{\nabla})[-1])=H^0(C, \textup{coker}(\bar{\nabla}))$, which is the length of the support of the cokernel. It must be zero by
\eqref{eqn_GrF1}, which implies the claim.
\end{proof}




\begin{corollary} (cf. \cite[Theorem 1.1]{EGW})
    The period map $P_{1,1}: F_{1,1}\to \bD/\Gamma$ from the moduli space of elliptic-elliptic surfaces is dominant.
\end{corollary}
\begin{proof}
    Since the domain and target of $P_{1,1}$ are both 10-dimensional, and the period map is unramified at a point $p\in F_{1,1}$, the image is a dense open subset of $\bD/\Gamma$.   
\end{proof}
As a consequence of Lemma \ref{lemma_IVHS}, the fibers of the period map have dimension at most one. See \cite{Ikeda,EGW} for some examples of positive dimensional fibers for $P_{1,1}$.

\section{Noether-Lefschetz Locus}\label{hurwitz}

In this section, we introduce the Hurwitz scheme underlying the family of elliptic-elliptic surfaces in Proposition \ref{prop_basechangeX'}, and identify the one-parameter subfamily of Proposition \ref{prop_rho=11family} as a component of the Noether-Lefschetz locus.
\subsection{Hurwitz scheme} Consider the Hurwitz moduli space $\mathcal H$ whose points are smooth curves $E$ of genus one, equipped with a degree 3 morphism $\pi:E \to \bP^1$ triply branched over $0,\infty\in \bP^1$ and simply branched over two other distinct points $r_1,r_2\in \bP^1$. We do not quotient by automorphisms of the target $\bP^1$.
\begin{proposition} The moduli functor of such triple covers $\pi:E\to \bP^1$ is representable by a connected \'{e}tale double cover of $\Sym^2(\C^*)\smallsetminus \Delta$, which is quasi-projective.
\end{proposition}
\begin{proof}
    There is a map from $\mathcal H \to \Sym^2(\C^*)\smallsetminus \Delta$ sending $[\pi:E \to \bP^1]$ to the divisor $r_1+r_2$ of simple branch points. By the Riemann Existence Theorem, the triple cover $\pi$ is uniquely determined by the locations of its branch points in $\bP^1$ and their monodromy data. More precisely, labeling the three sheets of a cover, we obtain a transitive monodromy action
    $$\rho:\pi_1(\C^*\smallsetminus \{r_1,r_2\},*) \to S_3,$$
    Choosing loops around each branch point $0,r_1,r_2,\infty$ based at $*$, their images under $\rho$ must be a 4-tuple of permutations
    $$(\sigma_0,\tau_1,\tau_2,\sigma_{\infty}),$$
    where $\sigma_0,\sigma_\infty$ are 3-cycles and $\tau_1,\tau_2$ are transpositions, and $\sigma_0\tau_1\tau_2 \sigma_\infty=1$. There are only two such 4-tuples, modulo simultaneous conjugation:
    $$((123),(12),(23),(123)),\,\, ((123),(12),(12),(132)).$$
    This implies that the map $\mathcal H \to \Sym^2(\C^*)\smallsetminus \Delta$ is \'{e}tale of degree two. Using the braiding action on 4-tuples, one checks that a small loop around $\Delta$ acts trivially on 4-tuples modulo simultaneous conjugation, and a small loop around $0+t$ or $\infty+t$ ($t\in \C^\times$) acts non-trivially on 4-tuples modulo simultaneous conjugation, so $\mathcal H$ is connected. The connected \'{e}tale double cover may be realized in the category of $\C$-schemes by taking the cyclic double cover $\epsilon$ of $\bP^2\simeq \Sym^2(\bP^1)$ branched at the union of two disinct lines $L_0+L_{\infty}$, and then removing the ramification locus and the pre-image of $\Delta$.
\end{proof}

There is a universal family $\mathcal E \to \mathcal H$ of genus one curves, equipped with a universal map $\pi: \mathcal E \to \bP^1$. The two parameter family of elliptic-elliptic surfaces from Proposition \ref{prop_basechangeX'} is given by the fibered product
$$\mathcal Y := \mathcal E \times_{\bP^1} X' \to \mathcal H.$$

\subsection{The subfamily $T$}

We have a diagram

     \begin{equation}
\begin{tikzcd}
 &  \mathcal{H}\arrow[d,"\epsilon"]\\
\C^*\smallsetminus\{1\}\arrow[r,"\phi"] \arrow[ur,dashed, "\tilde{\phi}"]& Sym^2(\C^*)\smallsetminus \Delta,
\end{tikzcd}
\end{equation}
where $\phi$ sends $a$ to the unordered pair of branch points $\{(1+\sqrt{a})^2,(1-\sqrt{a})^2\}=\{r_1,r_2\}$, or in the natural coordinates on $\Sym^2(\bP^1)\simeq \bP^2$,
$$\phi: a\mapsto (r_1r_2, r_1+r_2) = ((1-a)^2,2+2a)=(u,v)$$
Observe that $\Sym^2(\C^*)\cong \C^2\smallsetminus \{u=0\}$, and $\phi$ induces a map on fundamental groups
$$\phi_*:\pi_1(\C^*\smallsetminus\{1\})\to \pi_1(\Sym^2(\C^*))\cong \Z,$$
sending a loop around $1\in \C^*$ to twice the generator, and a loop around $0\in \C^*$ to the identity, so $\phi$ lifts to the Hurwitz scheme $\mathcal{H}$:
$$\tilde{\phi}:\C^*\smallsetminus\{1\}\to \mathcal{H}$$
The image of $\phi$ in $\Sym^2(\bP^1)\simeq \bP^2$ compactifies to a conic tangent to $L_0$ and $L_\infty$, so its pre-image under $\epsilon$ is isomorphic to two disjoint copies of $\C^*\smallsetminus \{1\}$. The lift corresponding to the subfamily $T$ of Proposition \ref{prop_rho=11family} is determined by
$$\lim_{a\to 0} \tilde{\phi}(a),$$
which by Proposition \ref{prop_LMHS_0} is located on the sheet labeled by the Hurwitz conjugacy class
$$((123),(12),(23),(123)).$$
The diagonal $\Delta\subset \Sym^2(\bP^1)$ is also a conic tangent to $L_0$ and $L_\infty$, so its pre-image in the double cover of $\bP^2$ contains two disjoint copies of $\C^*$. The Hurwitz space does not contain these curves by our original definition, but that definition can be extended to allow the simple ramification points to collide. On the first sheet, the limiting triple cover $E_0 \to \bP^1$ is Galois where $E_0$ is smooth with an automorphism of order 3. On the second sheet, the limiting triple cover $E'_0\to \bP^1$ is Galois where $E'_0$ is a nodal rational curve with an automorphism of order 3. The compactified Hurwitz scheme is no longer a fine moduli space because of the presence of these automorphisms.





\subsection{Proof of Main Theorem}
We now complete the proof of Theorem \ref{thm_NL}. 

\begin{theorem}
    There is a surface $S$ in the moduli space $F_{1,1}$ of elliptic-elliptic surfaces, a curve $T\subseteq S$ and a point $p\in T$, such that 
\begin{itemize}
    \item a general point of $S$ is an elliptic-elliptic surface with four singular fibers of type $I_3$, $\rho=10$, and Mordell-Weil rank $0$;
    \item a general point of $T$ is an elliptic-elliptic surface with four singular fibers of type $I_3$, $\rho=11$, and Mordell-Weil rank $1$;
     \item $p$ corresponds to the elliptic-elliptic surface from Theorem \ref{mainthm}, with two singular fibers of type $I_3$ and one of type $I_6$, $\rho=12$, and Mordell-Weil rank 1.
\end{itemize}
\end{theorem}

\begin{proof}
The family $\mathcal Y \to \mathcal H$ constructed above corresponds to a morphism $\xi:\mathcal H \to F_{1,1}$ whose image we call $S$. Let $T\subset S$ be the image of 
$$\xi\circ \tilde{\phi}:\C^* \smallsetminus \{1\} \to F_{1,1}.$$
Let $p = \xi\circ \tilde{\phi}(4)\in T$, which corresponds to the surface $Y_p$ of Proposition \ref{prop_maximalPicard1}. Since $Y_p$ has a different singular fiber configuration ($2I_3+I_6$) from the other elliptic surfaces ($4I_3$) in the one-parameter family over $\C^*\smallsetminus \{1\}$, the image $T$ must be a curve. Now $S$ must be irreducible, and its closure in $F_{1,1}$ contains $T$, as well as one $\C^*$ component of the diagonal passing through $[Y_0]$ of Proposition \ref{prop_maximalPicard2}. Therefore, $S$ is a surface.

By the infinitesimal Torelli theorem (Corollary \ref{cor_IVHS}) for all surfaces in the family $\mathcal Y \to \mathcal H$, $P_{1,1}(S)\subset \bD/\Gamma$ is a surface. Likewise, $P_{1,1}(T)$ must be a curve, so the generic Picard number is 11 and not higher; see Proposition \ref{prop_rho=11family}. Since $Y_p$ has maximal Picard number, $P_{1,1}(p)$ is a CM point. There are countably many other CM points in $P_{1,1}(T)$, and they are equidistributed with respect to the Hodge metric by \cite{tayou}.

\end{proof}

\begin{remark}\normalfont
In fact, one can show $\xi\circ \tilde{\phi}$ is generically 1-to-1 onto $T$. Since the degree of the $J$-map for the family of base curves $\{E_a\}$ is equal to 4 by equation \eqref{eqn_JEa}, the moduli map to $T$ has degree at most 4. For generic $j$, the 4 surfaces whose base curve has $J$-invariant $j$ have different values under $\phi$, so they map to different points in $F_{1,1}$. 
\end{remark}




\section{A Shioda Modular Surface}\label{sec_Shioda}
From the arithmetic point of view, the most special elliptic surfaces are those considered by Shioda in \cite{Shi72}. Let $\Gamma'\subset \SL_2(\Z)$ be a congruence subgroup which acts effectively on the upper half plane $\bH$, so in particular $-I\notin \Gamma$. The quotient $Y_{\Gamma'} =\bH/\Gamma'$ is a modular curve with finitely cusps. By adding in the orbits of $\Gamma$ acting on $\bP^1(\bbQ)$, we obtain a smooth proper curve $X_{\Gamma'}$ which compactifies $Y_{\Gamma'}$. Likewise,
$$(\C\times \bH) / (\Z^2\rtimes\Gamma') \to Y_{\Gamma'}$$
may be compactified to an elliptic surface $B_{\Gamma'} \to X_{\Gamma'}$, where the fibers over the cusps are singular and their Kodaira types are completely determined by the $\Gamma'$-stabilizer of the cusp in question.

Take $\Gamma'\subset \Gamma_0(11)$ to be the kernel of the quadratic character $(\Z/11\Z)^\times \to \mu_2$ given by
$$\chi \begin{pmatrix} a & b \\ c & d \end{pmatrix} := \left( \frac{a}{11} \right).$$
Then the modular curve $X_{\Gamma'}$ has two cusps, $i\infty$ and $0$, and $B_{\Gamma'}$ has singular fibers of types $I_1$ and $I_{11}$ over the respective cusps. Since $X_{\Gamma'}$ is a genus one curve, we obtain that $B_{\Gamma'} \in F_{1,1}$ because it has topological Euler charateristic 12, and it has a distinguished section provided by the origin in $\C$.

The fundamental line bundle is a particular line bundle of degree $L$ on $X_{\Gamma'}$ of degree 1. Its sections correpond to modular forms of weight 1 for $\Gamma'$. 

\begin{proposition}
    The canonical fiber $K_B$ of $B_{\Gamma'}$ is a smooth fiber of $B_{\Gamma'} \to X_{\Gamma'}$.
\end{proposition}
\begin{proof}
    Consider the integral binary quadratic form $Q$ given by
    $$Q(x,y) = x^2+xy+ 3y^2.$$
    Let $\theta(q)$ be the associated theta function
    $$\theta(q) = \sum_{n\geq 0} r_Q(n) q^n,$$
    where $r_Q(n)$ is the is the number of representations of $n$ as $Q(x,y)$. The Poisson summation formula implies that $\theta(q)$ is the Fourier expansion of a modular form with character $\chi$:
    $$\theta \in M_1(\Gamma_0(11), \chi).$$
    It suffices to locate the unique zero of $\theta$ in $X_{\Gamma'}$. To do so, we consider the Jacobi product cusp form
    $$h(q) = q \prod_{n\geq 1} (1-q^n)^2(1-q^{11n})^2 \in M_2(\Gamma_0(11)).$$
    
    The quotient $F(q)=\theta(q)^2/h(q)$ is a meromorphic modular form for $\Gamma_0(11)$ of weight 0. The Atkin-Lehner involution $w(\tau)= -1/11\tau$ acts on the modular curve $X_{\Gamma'}$, interchanging the two cusps, and by the Poisson summation formula, it has the property that
    \begin{equation}\label{al} F(w(\tau)) = F(\tau).\end{equation}
    Therefore, since $F$ has a simple pole at $q=0$ ($\tau =i\infty$), it also has a simple pole at $q=1$ ($\tau = 0)$. This means that $\theta(q)$ must be nonzero at both cusps. Furthermore, \ref{al} implies that the zero of $\theta$ must lie at a fixed point of $w$. There are 4 such fixed points in $X_\Gamma$, and a computation shows that $\tau = i/\sqrt{11}$ gives the zero.
\end{proof}
The $J$-map associated to $B_{\Gamma'}$ has degree 12, ramified only at the points $\tau=0,i,e^{2\pi i/3}$. The canonical fiber does not lie over any of these values, so the Kodaira-Spencer map $\kappa$ is nonzero there. 

\begin{corollary}\label{shiodatorelli}
    The infinitesimal Torelli theorem holds for $B_{\Gamma'}$.
\end{corollary}
\begin{proof}
    This follows immediately from Lemmas \ref{lemma_IVHS}, and \ref{lemma_KSextremal}, using the fact that $B_{\Gamma'}$ is an extremal elliptic surface.
\end{proof}

\section{Birational Involution}\label{sec_involution}
In this section, we study a birational involution on the extremal rational elliptic surface $X'$ with singular fibers of type $IV^*,I_3$, and $I_1$ (cf. Lemma \ref{lemma_X'}). It turns out that the trisections $E_a$ (cf. Proposition \ref{prop_trisection/Z3}) on $X'$ are the fibers of another (rational) projection. There is a natural rational involution $\tau$ on $X'$ that sends $E_a$ to a fiber of the standard projection. We extend $\tau$ to a regular involution on a certain blow-up of $X'$. This provides a model for the base curves of the three limiting surfaces ``$Y_0$'', ``$Y_1$'', and ``$Y_{\infty}$''. This will be used in Section \ref{sec_compactify} to compute the limiting mixed Hodge structure at the three boundary points.

\subsection{Extension of birational involution}

\begin{proposition}\label{prop_birinv}
    There is a birational involution $\tau:X'\to X'$ which permutes the two projections

 \begin{equation}\label{eqn_birinv}
\begin{tikzcd}
 & \arrow[loop above] X' \arrow[dl,dashed,"\pi_y"]\arrow[dr,"\pi_v"]\\
\mathbb P^1_y  && \mathbb P^1_v,
\end{tikzcd}
\end{equation}
where $\pi_v$ is the standard elliptic fibration, while the general fiber of $\pi_y$ is the trisection $E_y$.
\end{proposition}
\begin{proof}
 
We begin by finding an affine equation for the quotient of the Hesse pencil in the $t$-chart. In equation \eqref{eqn_number68}, we substitute $u=x^3$, $v=t^3$, and $w=xt$, to get the Weierstrass equation of the $\Z_3$-quotient:
\begin{equation}\label{eqn_affineX'}
    (y+v+3w-1)yv=w^3.
\end{equation}
The involution is given by $\tau:(y,v)\mapsto (v,y)$ and fixes the $w$-coordinate. Since the trisection $D_a$ on the Hesse pencil is defined by $y=a$ (cf. Definition \ref{def_Da}), its quotient trisection $E_a$ has an affine equation given by setting $y=a$ in \eqref{eqn_affineX'}, which is a general fiber of $\pi_y$.
\end{proof}

\begin{proposition} \label{prop_X'tilde}
The birational involution $\tau$ extends to a biregular involution $\tilde{\tau}$ on a smooth rational surface $\tilde{X'}$. In particular, $\tilde{X'}$ has elliptic fibrations in two ways. Its relatively minimal regular model with respect to either of the projections is isomorphic to the extremal rational elliptic surface $X'$ with singular fiber types $IV^*,I_3,I_1$. 
\begin{figure}[h!]
\centering
\includegraphics[width=0.7\textwidth]{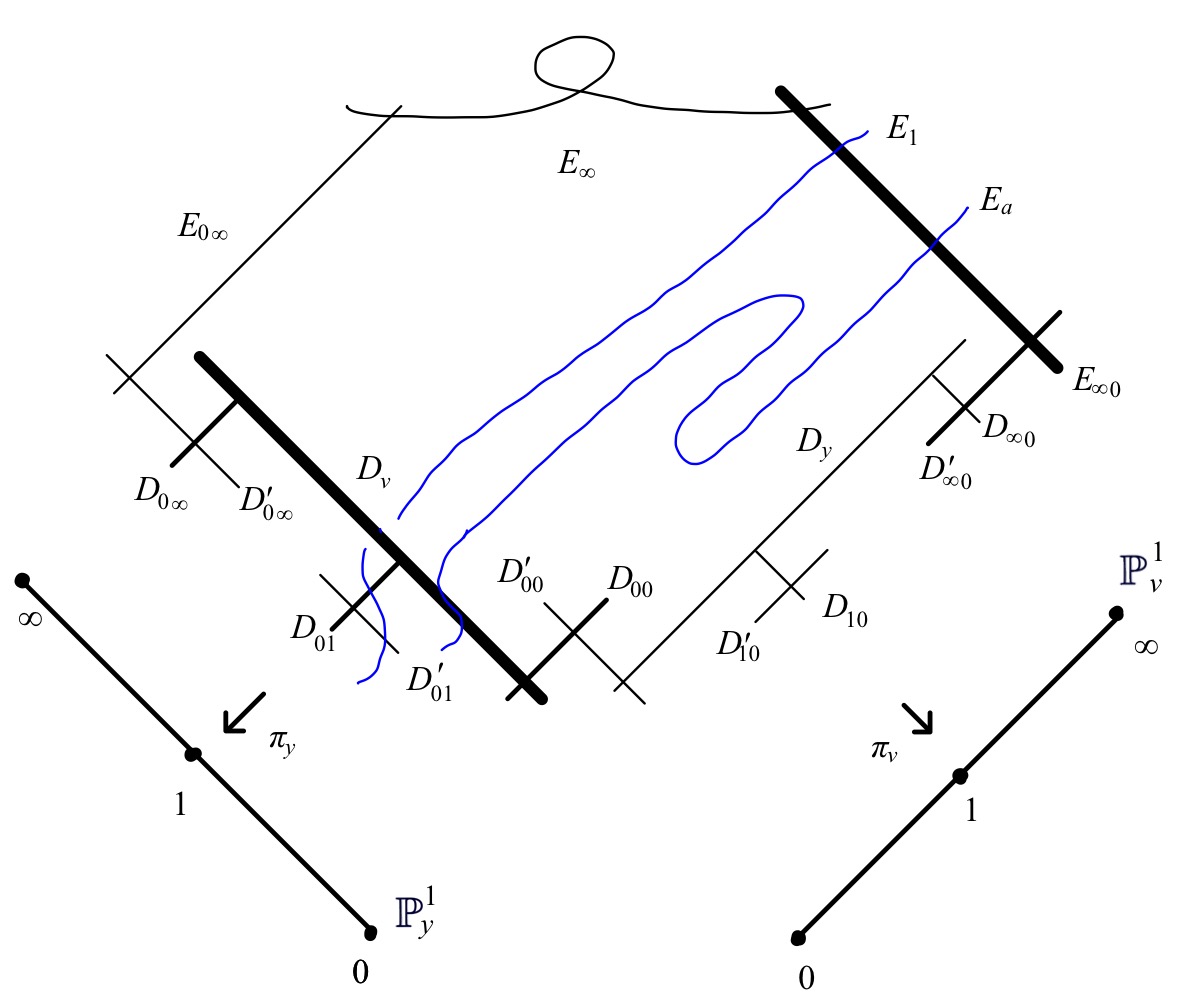}
\caption{Two elliptic fibrations on the blow-up of $X'$}
\label{figure_Xtildeprime}
\end{figure}

In the incidence diagram above, $\tilde{\tau}$ acts by the reflection about the vertical axis. The curves $E_{0\infty}$ and $D_y$ are sections of $\pi_v$. All the other black curves are supported on the fibers of $\pi_v$, with nonreduced components in bold. In particular,

\begin{itemize}
    \item $\pi_v^*(0)=3D_v+2(D_{00}+D_{01}+D_{0\infty})+D_{00}'+D_{01}'+D_{0\infty}'$ is the $IV^*$ fiber;
    \item $\pi_v^*(1)=D_{10}+D_{10}'+\tilde{\tau}(E_1)$ is the $I_3$ fiber;
    \item $\pi_v^*(\infty)=E_{\infty}+D_{\infty 0}+2D_{\infty 0}'+3E_{\infty 0}$ is a chain of rational curves, with $E_\infty$ the proper transform of the $I_{1}$ fiber on $X'$.
\end{itemize}

Here, $\pi_y^{-1}(a)=E_a$ is a general trisection (cf. Proposition \ref{prop_trisection/Z3}), meeting $D_v$ with intersection multiplicity 3 and $E_{\infty 0}$ transversely. The rational curve $E_1$ is also a trisection with respect to $\pi_v$ and intersects the $D_{01}$ and $D_{01}'$ transversely. Together, they form an $I_3$ fiber with respect to the projection $\pi_y$.

\end{proposition}
\begin{proof}
    One compactifies the affine variety \eqref{eqn_affineX'} in $\mathbb P^3$ with coordinates $[y:v:w:z]$ to get $$(y+v+3w-z)yv=w^3,$$ 
    
    which is a cubic surface $M$ with three $A_2$ singularities at $p_{0,0}=[0:0:0:1]$, $p_{1,0}=[1:0:0:1]$, and $p_{0,1}=[0:1:0:1]$. The three singularities are supported on the affine part \eqref{eqn_affineX'}. The curve at infinity given by $z=0$ is a nodal cubic curve, where the total space is smooth. Note that the equation is symmetric in $y$ and $v$, so the involution $\tau$ extends to $M$. We can resolve the three $A_2$ singularities, and the exceptional divisors $D_{00},D_{00}'$, $D_{01},D_{01}'$, and $D_{10},D_{10}'$ form three chains of $(-2)$-curves.

    Now, the projection $\pi_v$ (resp. $\pi_y$) has base locus supported at $b_v:=[1:0:0:0]$ (resp. $b_y:=[0:1:0:0]$) with length 3. To resolve the base locus of $\pi_v$, we can blow up the line $\{v=z=0\}$ in $\mathbb P^3$. Then $\pi_v$ extends to a morphism from the proper transform $M_v$, which has an $A_2$ singularity $p_{\infty,0}$ on the exceptional divisor $\bar{E}_{0\infty}$. Resolve the singularity. Then, the two exceptional curves $D_{0\infty}$ and $D_{0\infty}'$ are mapped to $v=0$, and the proper transform $E_{0\infty}$ of $E_{0\infty}$ is 1-to-1 onto $\mathbb P^1_v$ and becomes a section. We get a smooth surface $\tilde{M}$ with a morphism to $\mathbb P^1_v$. The fiber over $v=0$ consists of $D_v$ with multiplicity 3 as the proper transform of the line $\{v=w^3=0\}$ from the cubic surface $M$, together with three chain of rational curves (over $p_{1,0}, p_{0,0}$, $p_{\infty,0}$), which form a $IV^*$ configuration. The fiber over $v=1$ is of type $I_3$ (proper transform of the nodal cubic $\{v=1\}$ on $M$ together with exceptional curves over $p_{0,1}$). Finally, the fiber $E_{\infty}$ over $v=\infty$ comes from the nodal cubic curve $\{z=0\}$ on $M$. Therefore the surface is relatively minimal, and has singular fibers $IV^*, I_3, I_1$. By the classification of extremal rational elliptic surfaces, it is isomorphic to $X'$. One can resolve the base locus for $\pi_y$ symmetrically. Hence, the birational involution extends to a regular involution preserving the two projections.
\end{proof}

\begin{remark}\normalfont
    To get the elliptic fibration $X'\to\mathbb P^1_v$, one can blow down the three curves $E_{\infty 0}, D_{\infty 0}'$, and $D_{\infty 0}$ inductively. 
\end{remark}
\subsection{Degeneration of base curves}

\begin{corollary}  \label{cor_afamily}
    The one-parameter family of trisections $\{E_a\}$ for ${a\in \mathbb P^1\smallsetminus \{0,1,\infty\}}$ compactifies to the rational elliptic surface $\tilde{X}'$ with singular fibers $IV^*, I_3$ over $a=0,1$, and $I_1$ with three rational tails over $\infty$. Using $\pi_v$, the triple covers $E_a\to \mathbb P^1$ can be extended to the fibers over $a=0,1,\infty$. 
\end{corollary}
Note that the fibers over $a=0$ and $a=\infty$ are not reduced. We will replace those fibers with a semistable reduction with respect to the projection $\pi_v$. This will be discussed in some depth in Section \ref{subsec_a=0} and \ref{subsec_a=infty}. Our semistable reduction of the base curve is not relatively minimal at $a=\infty$; we need to keep a $(-1)$-curve in the central fiber to ensure that $\pi_v$ extends to a morphism.

\begin{corollary}
    The monodromy action on $H^1$ near $a=0$ has finite order, while the monodromy actions on $H^1$ near $a=1$ and $\infty$ have infinite order.
\end{corollary}
\begin{proof}
   This follows from the fact that $H^1(Y_a)\cong H^1(E_a)$ as weight one Hodge structures \cite[Lemma VII.1.1]{Miranda}, and the monodromy classification for Kodaira fibers \cite[p.36]{CoxZuc79}. 
\end{proof}

\section{Compactifying the Period Map}\label{sec_compactify}
In this section, we study the extension of the period map $P_{1,1}$ for the one-parameter family $T$ at the three limit points $a=0,1,\infty$, by working with local models for each degeneration. We consider first the family of the base curves $E_a$, together with their triple covers 
$$\pi:E_a \to \bP^1.$$

The monodromy action on $E_a$ near $a=0$ has finite order, so after a finite base change, there is a smooth semistable limit $E_0$. The surfaces $Y_a$ also have a smooth semistable limit there, which turns out to be the extremal elliptic-elliptic surface in Proposition \ref{prop_maximalPicard2}. Near $a=1$ and $a=\infty$, the monodromy action on $E_a$ has infinite order, so we will complete the family of surfaces with normal crossing limits. It turns out that at $a=1$, the limit has a component birational to a K3 surface, so the limiting mixed Hodge structure is pure; at $a=\infty$, all the strata of the limit are rational, so the limiting mixed Hodge structure has mixed weights. 


\subsection{Limit at $a=0$}\label{subsec_a=0}
In Corollary \ref{cor_afamily}, we found that the limiting elliptic curve at $a=0$ is singular of type $IV^*$. By \eqref{eqn_IV*basechange}, after taking a local 3-to-1 base change at the $a=0$ fiber, the central fiber can be replaced by a smooth elliptic curve with $J$-invariant zero. This suggests that we apply the same local base change of the family $E_a\times_{\mathbb P^1}X'$ of elliptic-elliptic surfaces (cf. Proposition \ref{prop_rho=11family}).
\begin{proposition}\label{prop_E0}
    After base change $b^3=a$, the cover $E_{b^3}\to \mathbb P^1$ limits to a Galois cover
    $$E_0\to \mathbb P^1,$$
totally branched at $0,1,\infty$, where $E_0$ is the genus one curve with $J$-invariant zero. 
\end{proposition}
\begin{proof}
Explicitly, by setting $\tilde{y}^3=y$ in \eqref{eqn_affineX'}, the surface becomes non-normal. By substitution $\tilde{w}=w/\tilde{y}$, we find the equation of the normalization:
$$(\tilde{y}^3+v+3\tilde{y}\tilde{w}-1)v=\tilde{w}^3.$$

The fiber at $\tilde{y}=0$ is the smooth elliptic curve $E_0:=\{v^2-v=\tilde{w}^3\}\cup\{\infty\}$, which has $J$-invariant zero. Moreover, the 3-to-1 map $E_{\tilde{y}^3}\to \mathbb P^1$ given by $\pi_v$ extends to the smooth limit at $\tilde{y}=0$:
$$E_0\to \mathbb P^1_v,\ (\tilde{w},v)\mapsto v.$$
It is totally branched at $v=0,1,\infty$, which are exactly the critical values of $\pi_v: X' \to \bP^1$.
\end{proof}

\begin{proposition}\label{prop_LMHS_0}
    The semistable limit of $\{Y_a\}$ for $a\to 0$ is a smooth extremal elliptic-elliptic surface with singular fibers of types $I_9$ and $I_3$ (cf. Proposition \ref{prop_maximalPicard2}). Consequently, the period map at $a=0$ extends to the interior of $\mathbb D/\Gamma$.
\end{proposition}
\begin{proof}
    By Proposition \ref{prop_E0}, the triple covers can be extended over the semistable limit. Therefore, the fiber product construction $E_{\tilde{y}^3}\times_{\mathbb P^1_v}X'$ extends to $\tilde{y}=0$, and the total space is equisingular (along the $IV^*$ fiber and nodal point on $I_1$ fiber) with respect to the $\tilde{y}$-direction. One can perform semistable reduction on these fibers as in Proposition \ref{prop_basechangeX'} and Section \ref{subsec_3-to-1IV*}, and get a smooth family $\{Y_{\tilde{y}}\}$, with $\tilde{y}\in \Delta_0$ a small disk centered at $0$. The limiting surface $Y_0$ is exactly the extremal elliptic-elliptic surface from Proposition \ref{prop_maximalPicard2}.
\end{proof}

\begin{remark}\normalfont
    One should note that the limiting surface $Y_0$ has no section of infinite order. This is because the limiting elliptic curve $E_0$ does not live in the rational elliptic surface $X'$ as a trisection. Rather, $E_0$ triple Galois cover of a (torsion) section of $X'$.

\end{remark}

From another point of view, the two simple branching points of $E_a\to \mathbb P^1$ are $(1\pm \sqrt{a})^2$, so when $a\to 0$, they collide to produce a total branching under the singular fiber $I_3$. This suggests that the semistable limit at $a=0$ is the elliptic surface in Proposition \ref{prop_maximalPicard2}.

\subsection{Limit at $a=1$.}
We will construct a normal crossings degeneration of $\{Y_a\}$ at $a=1$ after some birational modifications.

Note that Corollary \ref{cor_afamily} provides a model of the base of the limit ``$Y_1$''. Therefore, we can continuously extend the fiber product $E_a\times_{\mathbb P^1_v}X'$ to $a=1$. Since the semistable limit of the base curves consists of a triangle of rational curves, $D_{01}, D_{01}'$, and $E_1$, we should expect their fiber product $D_{01}\times_{\mathbb P^1_v}X'$, $D_{01}'\times_{\mathbb P^1_v}X'$, and $E_1\times_{\mathbb P^1_v}X'$ are components of ``$Y_1$'', at least birationally. We first show that the last component is birational to a K3 surface.



\begin{lemma}\label{lemma_K3}
    The relatively minimal regular model of the fibration $E_1\times_{\mathbb P^1_v}X'\to E_1$ is an elliptic K3 surface with Picard number of at least 19. It has singular fibers $4\times I_3$, $IV^*$, and $IV$, and a section of infinite order.
\end{lemma}
\begin{proof}
 The image of $E_1$ in the Weierstrass model \eqref{eqn_affineX'} of $X'$ is a nodal curve $\bar{E}_1$ with an affine equation
 \begin{equation} \label{eqn_E1}
    v^2+3wv=w^3.
 \end{equation}
 
 Moreover, the 3-to-1 map $E_1\to \mathbb P^1_v$ factors through the normalization $E_1\to \bar{E}_1$ and the projection of $\bar{E}_1$ onto the $v$-line. From the equation \eqref{eqn_E1}, $\bar{E}_1$ has two branches at $v=w=0$ — one branch is 2-to-1 onto the $v$-line, and the other one is 1-to-1 onto the $v$-line. In addition, $\bar{E}_1\to \mathbb P^1_v$ has a total branching at $v=\infty$ and another simple branching at $v=4$. Therefore, the composition $E_1\to \bar{E}_1 \to \mathbb P^1_v$ is 3-to-1, totally branched at $v=\infty$, and simply branched at $v=0$ and $4$.

 
 As in Proposition \ref{prop_basechangeX'}, the relatively minimal regular model $Z'$ of the fiber product $E_1\times_{\mathbb P^1_v}X'$ has four $I_3$ fibers, one $IV^*$ fiber, and one $IV$ fiber. The $IV$ fiber occurs over the simple ramification over $v=0$, and it arises from quadratic reduction of $IV^*$ fiber (cf. \eqref{eqn_IV*basechange}, Appendix \ref{subsec_2-to-1IV*}).   One checks that $Z'$ has topological Euler characteristic equal to 
 $$\chi=4\times \chi(I_3)+\chi(IV^*)+\chi(IV)=12+8+4=24,$$
 which implies that its fundamental line bundle has degree $2$. By the classification of elliptic surfaces with section \cite[Lemma III.4.6]{Miranda}, $Z'$ is a K3 surface.
\end{proof}



We will now construct a family $\{\tilde{Y}_a\}$, birational to the original family $\{Y_a\}$, near $a=1$, with normal crossing central fiber.

\begin{proposition}\label{prop_a=1ssmodel}
    There exists a smooth analytic 3-fold $\mathcal{Y}\to \Delta_1$ fibered over a disk $\Delta_1$ centered at $a=1$ such that:
    
    \begin{enumerate}
        \item The fiber over $a\neq 1$ is a smooth surface $\tilde{Y}_a$, whose relatively minimal regular model is the elliptic-elliptic surface $Y_a$;
        \item The fiber over $a=1$ is a normal crossing surface fibered over $I_3$;
        \item The fiber over $a=1$ has a smooth, reduced component $Z$ that is birational to an elliptic K3 surface.
    \end{enumerate}
    
\end{proposition}
\begin{proof}

We first construct a (singular) total space $W_1$ which contains the surface $\pi_y^{-1}(1)\times_{\mathbb P^1_v}X'$ as the central fiber. Then $\mathcal{Y}$ will be an appropriate desingularization of $W_1$.


Define $W_1$ to be the two-sided fiber product
 \begin{equation}\label{eqn_a=1_totalsing}
   W_1:=\Delta_1\times_{\mathbb P^1_y}\tilde{X'}\times_{\mathbb P^1_v}X'.  
 \end{equation}

Note that $\Delta_1\times_{\mathbb P^1_y}\tilde{X'}$ is a surface which models the base $\{E_a\}$ of the family $\{Y_a\}$ near $a=1$, and the second fiber product is to perform continuous base change of the rational elliptic surface $X'\to \mathbb P^1_v$. The projection $\pi_1:W_1\to \Delta_1$ has fiber over $a\neq 1$ equal to the (singular) elliptic surface $E_a\times_{\mathbb P^1_v}X'$ and birational to $Y_a$.
 
The fiber at $a=1$ consist of two constant families $D_{01}\times IV^*$, $D_{01}'\times IV^*$, and a smooth elliptic surface $Z\to E_1$ birational to a K3 surface (cf. Lemma \ref{lemma_K3}). They form a triangle as shown in the following picture.

\begin{figure}[h!]
\centering
\includegraphics[width=0.5\textwidth]{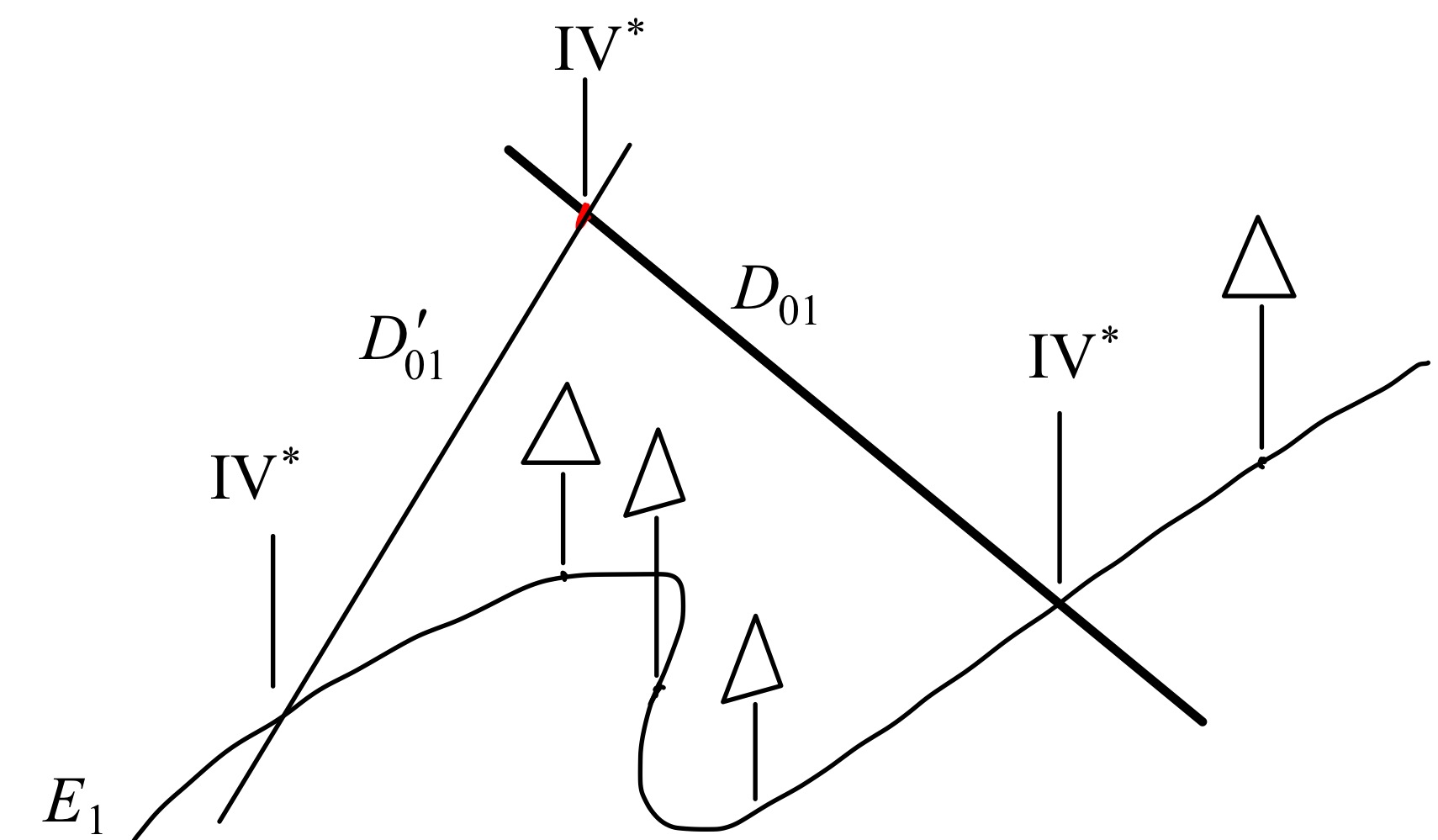}
\caption{Limiting model at $a=1$}
\label{figure_K3}
\end{figure}

Now, we need to desingularize $W_1$ relative to $\Delta_1$. Note $W_1$ is singular along the $IV^*$ fiber for each $a$ (where $W_1$ is non-normal), and is equisingular along the section marking the node of $I_1$ (where $W_1$ has transversal $A_2$ singularities). We first take the normalization $W_1'$ of $W_1$. Then by Appendix \ref{subsec_3-to-1IV*}, the general fiber of $W_1'\to \Delta_1$ has three terminal rational singularities (cf. Lemma \ref{lemma_twistedcubiccone}) and one $A_2$ singularity. They form four sections as $a$ varies in $\Delta_1$ and extends to $a=1$. Blowing up $W_1'$ along the four sections, we obtain a family $W_1''\to \Delta_1$ whose fiber at $a\neq 0$ is a smooth elliptic surface $\tilde{Y}_a$ birational to $Y_a$ (by blowing up nine smooth points of $Y_a$). Moreover, the singularity of $W_1''$ is supported on the $a=1$ fiber. Careful analysis reveals that the singularities are supported on the fiber over $D_{01}\cap D_{01}'$ and $D_v\cap D_{01}$. In particular, $W_1''$ is smooth at the component $Z$ which is birational to a K3 surface (cf. Lemma \ref{lemma_K3}). We choose a log resolution of $W_1''$ centered on the singular strata and obtain a smooth total space $\mathcal{Y}$, with the desired properties.
\end{proof}




\begin{remark}\normalfont
   Note that the central fiber of $\mathcal{Y} \to \Delta_1$ contains non-reduced components, which implies that the family in Proposition \ref{prop_a=1ssmodel} is not a semistable model. In fact, the monodromy on $H^2$ at $a=1$ has finite order and is non-trivial. Therefore, finding semistable reduction requires further base changes followed by a sequence of blow-ups and blow-downs. However, we do not pursue that in this paper.
\end{remark}

\begin{corollary}\label{cor_LMHS_1}
    The limiting mixed Hodge structure at $a=1$ is pure. The monodromy on $H^2$ at $a=1$ has finite order. As a consequence, the period image at $a=1$ is in the interior of $\mathbb D/\Gamma$.
\end{corollary}
\begin{proof}
Performing semistable reduction on the family $\mathcal{Y}\to \Delta_1$ produces a new family $\tilde{\mathcal{Y}}\to \tilde{\Delta}_1$. By Proposition \ref{prop_a=1ssmodel}, there is a component of the central fiber isomorphic to $Z$, which supports a K3-type Hodge structure. The Clemens-Schmid sequence \cite{ClemensSchmid} reads:
\begin{equation}\label{eqn_ClemensSchmid}
    H_4(\tilde{\mathcal{Y}}_1)\to H^2(\tilde{\mathcal{Y}}_1)\xrightarrow{\alpha} H^2_{lim}\xrightarrow{N} H^2_{lim},
\end{equation}
where $\tilde{\mathcal{Y}}_1$ is the fiber at $a=1$, and $N$ is the logarithm of the local monodromy operator on $H^2$ of the nearby fiber. 
Since $H_4(\tilde{\mathcal{Y}}_1)$ is Hodge-Tate, the composition
$$H^2(Z)\to  H^2(\tilde{\mathcal{Y}}_1)\to H^2_{lim}$$ 
is injective on $H^{2,0}$ and $H^{0,2}$. In particular, the image of $N$ does not shift the weight, so the monodromy action on $H^2$ near $a=1$ from Proposition \ref{prop_a=1ssmodel} has finite order.
\end{proof}

    

\subsection{Limit at $a=\infty$.}\label{subsec_a=infty}
To compute the limiting mixed Hodge structure at $a=\infty$, we will find a semistable model. We first extend the triple covers $E_a\to \mathbb P^1_v$ to $a=\infty$, and then extend the base change construction to $a=\infty$, resulting in a degeneration with singular total space. Then we resolve those singularities. 

By Corollary \ref{cor_afamily}, the limiting base curve $\pi_y^{-1}(\infty)$ is not reduced. We will first find a reduced model of $\{E_a\}$ at $a=\infty$. This amounts to finding a semistable reduction, except that we want to keep a $(-1)$-curve in the central fiber to ensure that the 3-to-1 map $\pi_v$ is still a morphism.

\begin{proposition}\label{prop_a=infty_base}
    Let $\tilde{\Delta}_{\infty}\to \Delta_{\infty}$ be the 3-to-1 map of a disk totally branched at $\infty$, then there is a smooth surface $\mathcal{E}$ together with projection $\mathcal{E}\to \tilde{\Delta}_{\infty}$ such that 
    \begin{enumerate}
        \item any fiber over $b\in \tilde{\Delta}_{\infty}\smallsetminus\{\infty\}$  is a smooth elliptic curve $E_b$
        \item over $b=\infty$, the fiber is reduced and consists of a smooth rational curve $R$ meeting an $I_3$ configuration transversely;
        \item The 3-to-1 map $E_a\to \mathbb P^1_v$ extends to $R\cup I_3$.
    \end{enumerate}
\end{proposition}
\begin{proof}
    By Proposition \ref{prop_X'tilde}, the fiber of $\pi_y$ at $\infty$ has multiplicity 3 along the component $E_{0\infty}$. After taking a 3-to-1 base change $\tilde{\Delta}_{\infty}\times_{\Delta_{\infty}}\tilde{X}'$
and normalizing, the central fiber becomes $E_{\infty}+R+2\tilde{D}_{0\infty}'+\tilde{D}_{0\infty}$, where $R\to E_{0\infty}$ is 3-to-1 and totally branched at two points, so $R$ is a rational curve. The total space has an $A_2$ singularity at the node of $E_{\infty}$ and an isolated rational singularity at $\tilde{D}_{0\infty}\cap\tilde{D}_{0\infty}'$, as in Lemma \ref{lemma_twistedcubiccone}. The minimal resolution of this $A_2$ singularity produces an $I_3$ fiber. For the other singularity, the reduction process is described in Appendix \ref{subsec_3-to-1IV*} — we blow up the singularity, contract the proper transforms of $\tilde{D}_{0\infty}$, $\tilde{D}_{0\infty}'$, and the exceptional curve to get a smooth total space $\mathcal{E}$ with reduced central fiber consisting of a smooth rational curve and a nodal curve meeting transversely.   

For the last statement, the morphism $\pi_v$ descends to the blow down $\tilde{X}'\to \mathscr{X}$ which contracts two $(-2)$-curves $D_{0\infty}$ and $D_{0\infty}'$. By Zariski's Main Theorem, $\mathcal{E}$ can be constructed by base change $\tilde{\Delta}_{\infty}\times_{\Delta_{\infty}}\mathscr{X}$ and normalization, so the morphism $\pi_v$ extends to $\mathcal{E}$.
\end{proof}

The central fiber $R+I_3$ of the family $\mathcal{E}\to \tilde{\Delta}_{\infty}$ at $\infty$ will be the base curve of the semistable limit ``$Y_{\infty}$''. 


\begin{proposition}\label{prop_a=infty_ssmodel}
    There exists a smooth analytic 3-fold $\mathcal{Y}^{\infty}$ together with a holomorphic morphism $\psi:\mathcal{Y}^{\infty}\to \tilde{\Delta}_{\infty}$ such that
      \begin{enumerate}
        \item The fiber over $b\neq \infty$ is the elliptic-elliptic surface $Y_b$;
        \item The fiber over $b=\infty$ is reduced, simple normal crossing, and each component is a rational surface;
        \item The fiber over $b=\infty$ has dual complex is homotopy equivalent to a pinched $2$-torus. 
        \end{enumerate}
\end{proposition}
\begin{proof}


 Using the model $\mathcal{E}$ constructed in Proposition \ref{prop_a=infty_base}, we take the fiber product
$$W_{\infty}:=\mathcal{E}\times_{\mathbb P^1_v}X'.$$
    The fibers of the natural projection  $\pi_{\infty}:W_{\infty}\to \tilde{\Delta}_{\infty}$ over $b\neq \infty$ are (singular) elliptic surfaces $E_a\times_{\mathbb P^1_v}X'$ with $IV^*, 3\times I_3$, and $I_1$ fiber and birational to $Y_a$. The fiber of $\pi_{\infty}$ at $\infty$ consists of four components arising from the pullback $X'$ to the four rational curves $R+I_3$. One component fibered over $R$ is a (singular) rational elliptic surface with singular fibers $IV^*,3\times I_3$, and $I_1$. The other three form a constant nodal curve fibration $I_1\times I_3$.


The total space $W_{\infty}$ is non-normal along the $IV^*$-locus for each $b\in \tilde{\Delta}_{\infty}$. In addition, it has a transversal $A_2$ singularity along the node of the $I_1$ fiber for each of $b\in \tilde{\Delta}_{\infty}$.

We first perform birational modifications to realize condition (1). Note that the 3-to-1 map $R\to \mathbb P^1_v$ is totally branched at $v=0$ and locally deforms to the nearby $E_a\to  \mathbb P^1_v$. Hence, the total space $W_{\infty}$ is equisingular along the  $IV^*$-locus. We can perform reduction uniformly as in Figure \ref{figure_IV*basechange}. In the end, the $IV^*$-locus is replaced by a smooth family of elliptic curves for each $b\in \tilde{\Delta}_{\infty}$. 

Now that the total space $W_{\infty}'$ is normal, we blow up the section which marks the transversal $A_2$ singularity, to arrive at a family  $W_{\infty}''\to \tilde{\Delta}_{\infty}$, whose fiber over $b\neq \infty$ is the smooth elliptic elliptic surface $Y_b$. However the total space $W''_\infty$ still has singularities over $\infty$.

     

\begin{figure}[h!]
\centering
\includegraphics[width=0.6\textwidth]{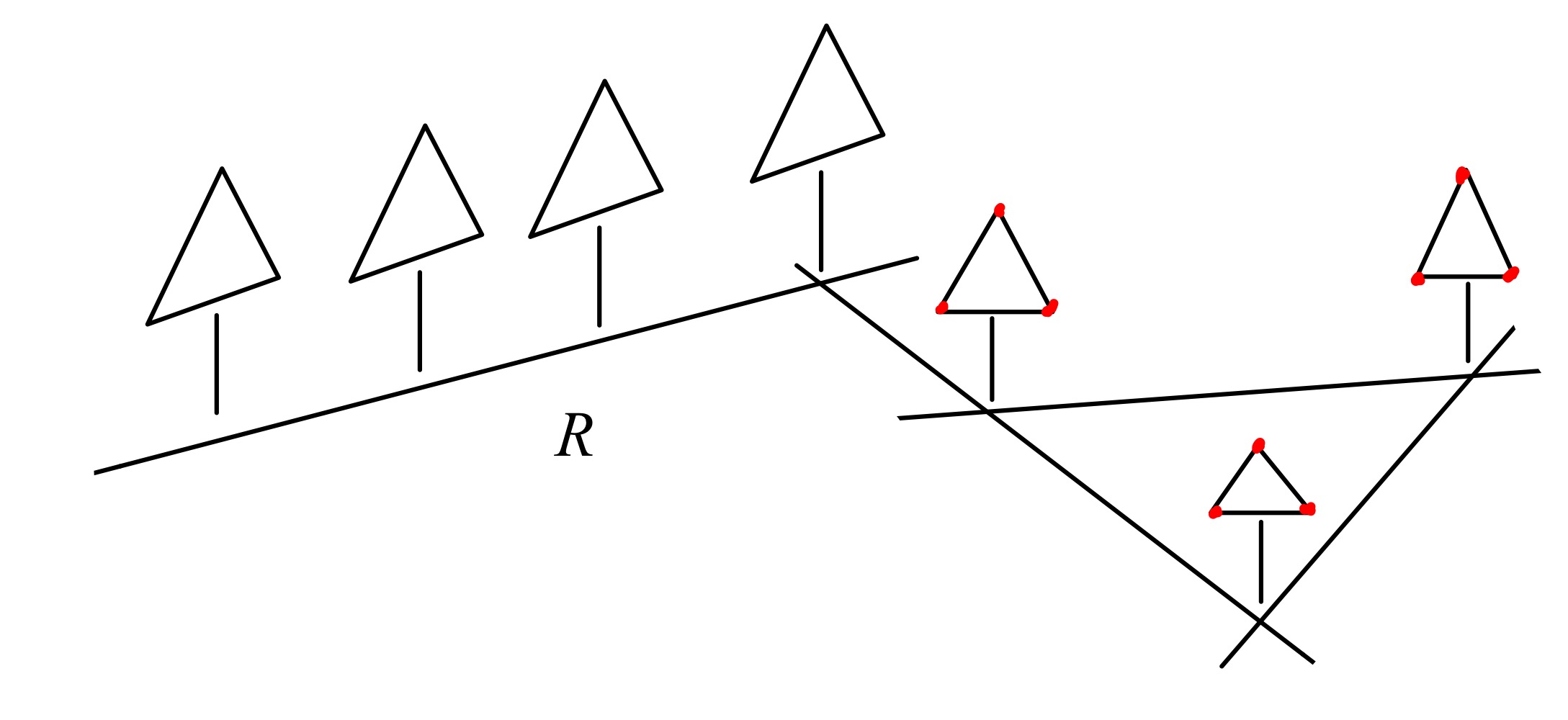}
\caption{Limiting model at $a=\infty$}
\label{figure_Limit_infty}
\end{figure}


The central fiber has four components, pictured above. The first is isomorphic to the Hesse pencil over $R$ (cf. Remark \ref{rmk_P1-to-P1}). The other three components given by $I_1\times I_3$. To resolve the singularities in the total space, we blow up the nodal locus on the constant $I_1$ fibration. This gives $I_3\times I_3$, and the total space has isolated singularities with analytic equation $xy=zw$ located at the product of two nodes. Blow up these singularities (there are 9 of them shown in Figure \ref{figure_Limit_infty}), and the total space becomes smooth, and the central fiber is reduced and simple normal crossing.

The dual complex of the $I_3\times I_3$ components is a 2-torus $\mathcal{T}^2$. The rational surface over $R$ adds a cone over a non-trivial loop on $\mathcal{T}^2$. The 9 new exceptional components lead to a triangular subdivision of $\mathcal{T}^2$. Therefore, the homotopy type of the dual complex of the central fiber is a pinched 2-torus.




\end{proof}

\begin{proposition}\label{prop_LMHS_infty}
    The limiting mixed Hodge structure at $a=\infty$ has type III. In particular, Hence, the monodromy action on $H^2$ at $\infty$ has infinite order.
\end{proposition}
\begin{proof}


From Proposition \ref{prop_a=infty_ssmodel}, we have a semistable model for the family $\{Y_{a}\}$, after taking a 3-to-1 base change. Since the central fiber $Y_{\infty}$ consists of only rational surfaces and meeting at rational curves, the mixed Hodge structure $H^2(Y_{\infty})$ is Hodge-Tate and is only supported on $(1,1)$ and $(0,0)$. We have a Clemens-Schmid sequence similar to \eqref{eqn_ClemensSchmid}, except that the image of $\alpha$ contains a nontrivial $(0,0)$ part. Therefore, the limiting mixed Hodge structure is type III.



\begin{align*}
\begin{tikzpicture}[scale=0.7]
    \draw[-,line width=1.0pt] (0,0) -- (0,3);
    \draw[-,line width=1.0pt] (0,0) -- (3,0);
\fill (2,2) circle (3pt);
\fill (1,1) circle (3pt);
\fill (1.2,0.8) circle (3pt);
\fill (0,0) circle (3pt);
\draw [-stealth](1.8, 1.8) -- (1.2,1.2) [line width = 1.0pt];
\draw [-stealth](0.8, 0.8) -- (0.2,0.2) [line width = 1.0pt];
\node at (1.7,-1) {LMHS at $a=\infty$};
\end{tikzpicture}
\end{align*}

In particular, the monodromy action has infinite order.
\end{proof}








\appendix

\section{Semistable Reduction of $IV^*$ fiber}\label{sec_app}
In this Appendix, we will work out semistable reduction for type $IV^*$ singular fibers under a base change of degree 2 and 3, respectively, and verify \eqref{eqn_IV*basechange}. As an application, this provides an alternative way to compute the degree of fundamental line bundle of the elliptic-elliptic surface arising from base change of a rational elliptic surface with an $IV^*$ fiber (cf. Proposition \ref{prop_basechangeX'}).  

\subsection{Quadratic base change}\label{subsec_2-to-1IV*}
An $IV^*$ fiber consists of 7 rational components: $E$ has multiplicity three, $D_1, D_2$ and $D_3$ have multiplicity two, while $C_1, C_2$, and $C_3$ have multiplicity one.

According to the principle of semistable reduction \cite[p.125]{HM98}, the 2-to-1 base change plus normalization has the effect of taking a branched cover of the components of odd multiplicity. Since the branch divisor is smooth, the resulting surface is smooth, and the pre-images of $C_i$ and $E$ have self-intersection $-1$. The components $D_i$ of multiplicity two undergo a 2-to-1 cover branched at the two intersection points with the $C_i$ and $E$. Lastly, we contract the $(-1)$-curves to obtain the fiber of type $IV$.

\begin{figure}[h!]
\centering
\includegraphics[width=0.7\textwidth]{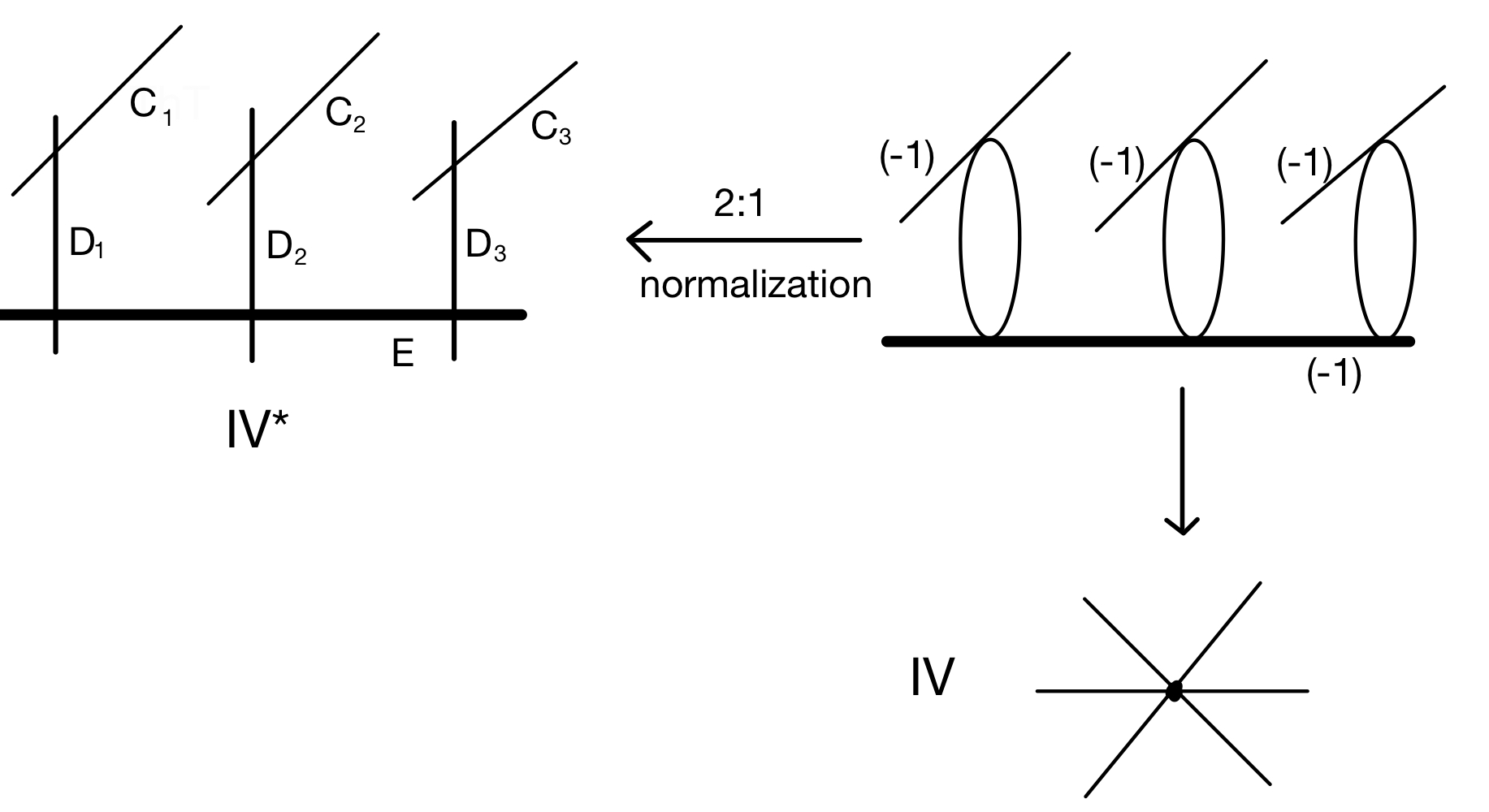}
\caption{2-to-1 reduction of $IV^*$}
\label{figure_IV*basechange2-to-1}
\end{figure}

\subsection{Cubic base change} \label{subsec_3-to-1IV*}
In the same way, according to the principle of semistable reduction \cite[p.125]{HM98}, the 3-to-1 base change plus normalization has the effect of taking a branched cover of the components with multiplicity one or two. The multiplicity three fiber undergoes a 3:1 cover totally branched at 3 points, so $g(E')=1$, and $E'\cdot E'=-6$.

\begin{figure}[h!]
\centering
\includegraphics[width=0.7\textwidth]{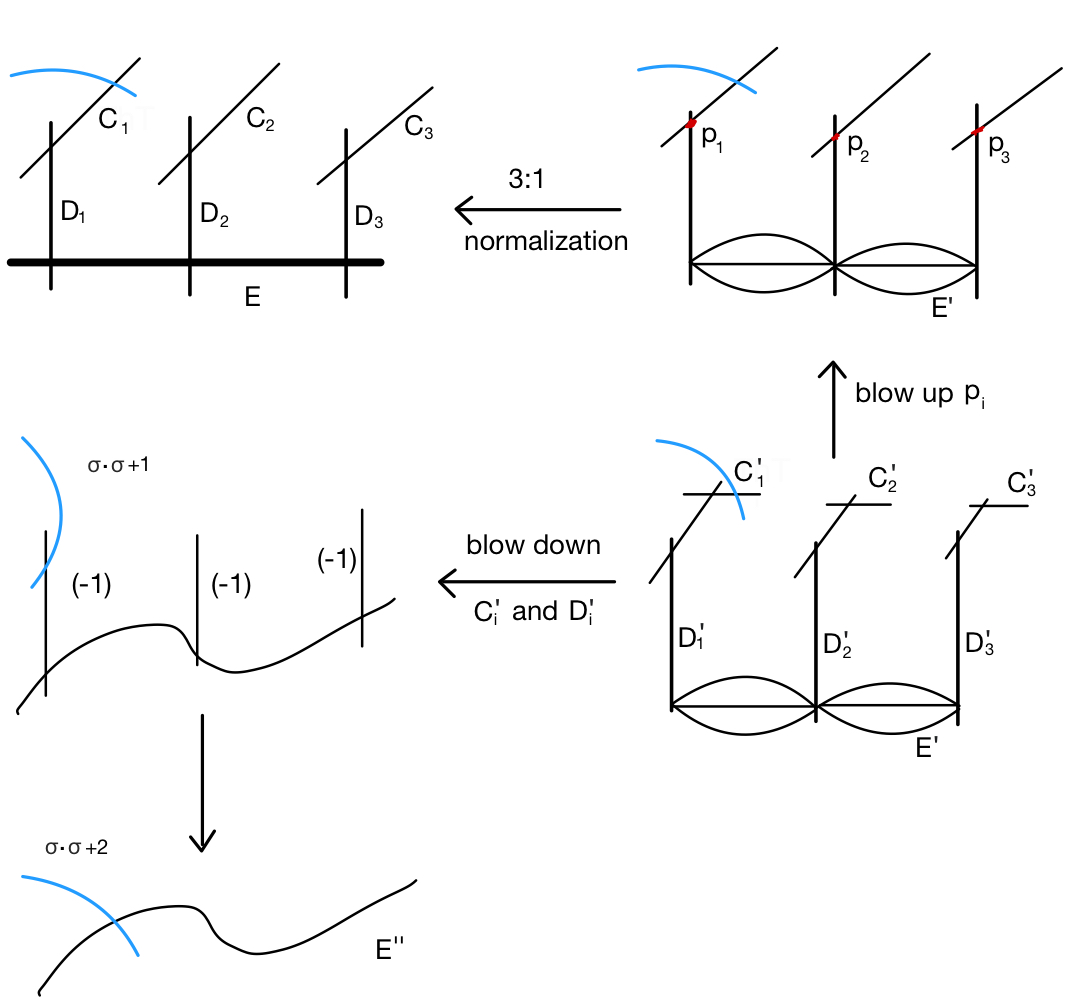}
\caption{3-to-1 reduction of $IV^*$}
\label{figure_IV*basechange}
\end{figure}

Next, the situation becomes more subtle as the total space is still not smooth — it has an isolated singularity at the intersection points of the multiplicity one and multiplicity two components. 

\begin{lemma} \label{lemma_twistedcubiccone}
 The singularity $p_i$ is isomorphic to the affine cone of twisted cubic.
\end{lemma}
\begin{proof}
    This amounts to studying the normalization of the ring $R=\C[x,y,t]/(t^3-x^2y)$. One observes $u^2=yt$ is an integral relation for $u=xy/t$, so we have $ut=xy$. Now multiply $t^3=x^2y$ by $y$ on both sides, and we get $u^2=yt$. Divide $t^3=x^2y$ by $t$ on both sides, we get $t^2=xu$. Therefore, we have three relations $ut=xy,u^2=yt,t^2=xu$ in the ring $\C[x,y,t,u]$. Since these relations define the ideal of the affine cone of twisted cubic, which is already normal, they generate the ideal of the normalization of $R$.
\end{proof}

Consequently, we blow up $p_i$ and get a $(-3)$-curve $E_i$. To compute the self-intersections of $C_i'$ and $D_i'$, we use the projection formula. Let $\pi$ be the composition of blow up, normalization, and base change. Then $\pi^*D_i=3D_i'+E_i$, $\pi^*C_i=3C_i'+E_i$. Then $-2=D_i\cdot D_i=\pi^*D_i\cdot D_i'=(3D_i'+E_i)\cdot D_i=3D_i'^2+1$, so $D_i'^2=-1$, and similar for $C_i'$. Now we blow down the $(-1)$-curves $C_i'$ and $D_i'$. The image of $E_i$ has self-intersection -1. Blow them down, leaving an elliptic curve $E''$ with self-intersection 0, which is the $I_0$ fiber. 


\begin{proposition}
    The fundamental line bundle $\mathcal{L}$ of the elliptic-elliptic surface $Y\to E$ that arises from the base change of rational elliptic surface with an $IV^*$ fiber in Proposition \ref{prop_basechangeX'} has degree $\deg(\mathcal{L})=1$.
\end{proposition}
\begin{proof}
    One can keep track of the change in the degree of the fundamental line bundle; it is the $-1$ times the self-intersection of a section $\sigma'$ \cite[Corollary 5.45]{SchShi19}. The section $\sigma'$ on $X'$ has self-intersection $\sigma'\cdot \sigma'=-1$ and passes through one of the multiplicity one components of the $IV^*$ fiber, say $C_1$. Its base change $\sigma:=\sigma'\times_{\mathbb P^1}E$ is away from the singular locus of the fiber product $X'\times_{\mathbb P^1}E$, its self-intersection is unaffected by the normalization, and $\sigma\cdot \sigma=-3$ by the projection formula. One can track the self-intersection of the section as we go through the reduction process. It only changes when we blow down curves intersecting $\sigma$. This happens twice, so the self-intersection of the image section $\sigma_Y$ on $Y$ is $\sigma_Y^2=\sigma^2+2=-1$, which implies that the fundamental line bundle has degree $1$.  
\end{proof}



\section{Hesse configuration}\label{sec_Hesse}
In this section, we want to find the projection of the trisection $D_a$ to $\mathbb P^2$ with Hesse configuration consisting of 9 points and 12 lines (cf. \cite{ArtDol}). Each point lies on 4 lines and each line contains 3 points. The 9 points arise from contracting the 9 (three-torsion) sections, and the 12 lines arise from the components on the four $I_3$ fibers.

We first find the curve class of $D_a$ in the Hesse pencil when $a\in \C^*\smallsetminus \{1\}$.

\begin{proposition}
    Let $X$ be the Hesse pencil. The trisection $D_a$ is linearly equivalent to $3C_0+3F$, where $C_0$ is the zero section and $F$ is a general fiber.
\end{proposition}
\begin{proof}
    By Proposition \ref{prop_Sa_zerosection}, $D_a$ has arithmetic genus 4. The adjunction formula gives $2g_a(D_a)-2=D_a\cdot D_a+K_X\cdot D_a$, and using the fact that $D_a$ is trisection, we get $D_a\cdot D_a=9$.

    On the other hand, $X$ is rational, so rational equivalence is the same as numerical equivalence. By Shioda-Tate, $H^2(X,\Z)$ has a basis given by $C_0,F$, and eight classes that are exceptional curves resolving the four $A_2$ singularities over $1,\zeta_3,\zeta_3^2$, and $\infty$ in the Weierstrass equation \eqref{eqn_number68}. Since by Proposition \ref{prop_Sa_zerosection}, the trisection $D_a$ is disjoint from these singularities, and both $C_0$ and $F$ are disjoint from these exceptional curve, so we can express
    $$[D_a]=\alpha [C_0]+\beta[F].$$

    Using $[D_a]\cdot [F]=3$ and $[D_a]\cdot [C_0]=0$, we have $\alpha=\beta=3$, so $[D_a]=3 [C_0]+3[F]$.
\end{proof}

Let $C_i$ for $1\le i\le 8$ be the nonzero 3-torsion sections. Let $H_a$ be the image of trisection $D_a$ in $\mathbb P^2$.

\begin{corollary}
  $D_a\cdot C_{i}=3$ for $1\le i\le 8$. 
\end{corollary}
\begin{proof}
    This just follows from that $C_{i}\cdot D_a=0$ and $C_{i}\cdot F=1$.
\end{proof}

\begin{proposition}
The image of $D_a$ in the $\mathbb P^2$ is a curve of degree $9$ and passes through 8 out of 9 points in the Hesse configuration and is singular there. Seven of the singularities are at least ordinary triple points, and the last singularity has local equation analytically equivalent to $x^3+y^6=0$.
\end{proposition}

\begin{proof}
Since $X$ is obtained from blowing up 9 points on $\mathbb P^2$, $K_X=-3\pi^*L+C_0+\cdots+C_8$, where $L$ is a general line on $\mathbb P^2$. We know $K_X\cdot D_a=-3$ since it has degree $3$. In the previous computation, we found that $C_0\cdot D_a=0$ and $C_i\cdot D_a=3$ for all $1\le i\le 8$, so
$$-3=-3\pi^*L\cdot D_a+3\times 8,$$
where $L$ is a line on $\mathbb P^2$ away from the blown up points. By the projection formula, 
$$L\cdot \pi_*(D_a)=\pi^*L\cdot D_a=9.$$
Hence, the plane curve $H_a$ has degree 9.

Among the nonzero 3-torsion sections, there is only one that meets the trisection $D_a$ at its ordinary triple point singularity, by the computation in Proposition \ref{prop_Sa_zerosection}. The singularity corresponding to the contraction of the ordinary triple point singularity of $D_a$ has local equation $x^3+y^6=0$. The other 3-torsion sections intersect $D_a$ at three distinct points (for general $a$). Hence, the contraction to $\bP^2$ produces an ordinary triple point.
\end{proof}

\begin{remark} \normalfont
The arithmetic genus of $H_a$ is $g_a(H_a)=\frac{8\cdot 7}{2}=28$, so its geometric genus is
$g(H_a)=28-7\delta+\delta'$
where $\delta$ and $\delta'$ are the delta invariants of the ordinary triple point and $x^3+y^6=0$, respectively. These two singularities have Milnor numbers $4$ and $10$, and they each have three branches. By the Milnor-Jung formula, $\delta=3$ and $\delta'=6$, therefore the $g(H_a)=1$. This agrees with Proposition \ref{prop_trisection} that there is a composition 
$$\tilde{D}_a\to D_a\to H_a$$
where the desingularization $\tilde{D}_a$ is a smooth genus one curve.
\end{remark}




\bibliographystyle{alpha}
\bibliography{bibfile}

\end{document}